# On a Fast, Robust Estimator of the Mode
## Comparisons to Other Robust Estimators with Applications

David R. Bickel [1]

Medical College of Georgia; Office of Biostatistics and Bioinformatics; USA
Email: bickel@prueba.info
URL: www.davidbickel.com

Rudolf Frühwirth

Austrian Academy of Sciences; Institute of High Energy Physics
Nikolsdorfer Gasse 18; A-1050 Wien; AUSTRIA
Email: fru@hephy.oeaw.ac.at
URL: wwwhephy.oeaw.ac.at/u3w/f/fru/www/

**Abstract**: Advances in computing power enable more widespread use of the mode, which is a natural measure of central tendency since, as the most probable value, it is not influenced by the tails in the distribution. The properties of the half-sample mode, which is a simple and fast estimator of the mode of a continuous distribution, are studied. The half-sample mode is less sensitive to outliers than most other estimators of location, including many other low-bias estimators of the mode. Its breakdown point is one half, equal to that of the median. However, because of its finite rejection point, the half-sample mode is much less sensitive to outliers that are all either greater or less than the other values of the sample. This is confirmed by applying the mode estimator and the median to samples drawn from normal, lognormal, and Pareto distributions contaminated by outliers. It is also shown that the half-sample mode, in combination with a robust scale estimator, is a highly robust starting point for iterative robust location estimators such as Huber's M-estimator. The half-sample mode can easily be generalized to modal intervals containing more or less than half of the sample. An application of such an estimator to the finding of collision points in high-energy proton-proton interactions is presented.

---

[1] Current address: Pioneer Hi-Bred International, Bioinformatics & Exploratory Research, 7250 NW 62nd Ave., PO Box 552, Johnston, Iowa 50131-0552

# 1 Introduction

An estimator is considered robust if it can be applied to samples drawn from a large class of distributions and if it is insensitive to outliers. The last three decades have seen a proliferation of research on robust descriptive statistics. This increased interest is due to more access to powerful computers and a growing recognition of the need to analyze data when the best model is unknown or when the data may be corrupted with outliers. In particular, the robustness of estimators of various measures of location has received much attention.

Staudte and Sheather [1] call the estimand $\mu(X)$ a *measure of location* if:

$$
\begin{aligned}
&\text{(i)} \quad \forall a > 0, b : \mu(aX + b) = a\mu(X) + b \\
&\text{(ii)} \quad \mu(-X) = -\mu(X) \\
&\text{(iii)} \quad X \geq 0 \implies \mu(X) \geq 0,
\end{aligned}
\tag{1}
$$

where $X$ is a random variable of a family of continuous distributions. Estimands satisfying these conditions include the mean, trimmed mean, median, linear combination of two quantiles, and mode. Bickel and Lehmann [2] would also require that

$$\text{(iv)} \quad \mu(X) \geq \mu(Y)$$

if $X$ is stochastically larger than $Y$, but that would rule out the mode as a measure of location [1]. Thus, condition (iv) is too restrictive for many applications:

1. the mode is more of a measure of central tendency than are measures that always satisfy (iv) since a difference only in the tails of distributions can lead to stochastic ordering but not to a difference in the mode;

2. although the mode does not satisfy (iv) in general, the mode does satisfy it for most two-tailed distribution families, e.g., the symmetric distribution families and the lognormal family;

3. the mode often satisfies (iv) for observed data with observations on both sides of the mode;

4. regardless of whether (iv) is satisfied, the mode has intuitive appeal as a measure of location since it is the value of $X$ with the maximum probability.

For these reasons, the definition of Eq. (1) is used herein to enable the study of the mode as a measure of location.

Because the median is the measure of location most often chosen when a robust measure is needed, it will be compared to the mode throughout this paper. Bickel and Lehmann [2] point out that since measures of location for a distribution differ only by a shift in the origin, they can be compared by comparing their estimators. As will be seen, a robust estimator of the mode performs better than the sample median in a variety of situations.

Many estimators of the mode have been suggested. For example, to estimate the mode $M$ of a unimodal distribution from the sample $\{x_1, x_2, \ldots, x_n\}$, with $x_1 \leq x_2 \leq \ldots \leq x_n$, Grenander [3] proposed the direct estimator

$$M^*_{p,k} = \frac{\frac{1}{2}\sum_{i=1}^{n-k}(x_{i+k}+x_i)/(x_{i+k}-x_i)^p}{\sum_{i=1}^{n-k} 1/(x_{i+k}-x_i)^p}, \quad 1 < p < k. \qquad (2)$$

Hall [4] proved the asymptotic normality of $M^*_{p,k}$ for $k > 2p$. A drawback of estimators of this type is that they are highly influenced by outliers. Section 2 features an estimator that overcomes these limitations; other robust estimators are also described. Various estimators of the mode are compared in Section 3. The featured estimator is applied to observed data in Section 4 and conclusions appear in Section 5.

## 2  Robust estimators of the mode

### 2.1  Half-sample mode

Consider a continuous distribution $F$ that only has one mode, $M$, and let $(x_i)_{i=1}^n$ be a vector of $n$ random numbers drawn from that distribution, ordered such that $x_1 \leq x_2 \leq \ldots \leq x_n$. Then the following procedure describes how to compute the *half-sample mode* (HSM) $\widehat{M}$, an estimate of $M$, when the sample size is an integer power of two ($n = 2^m$). (The algorithm is generalized in the Appendix to allow $n$ to be any positive integer.)

1. First, find the smallest interval that contains $n/2$ points from the sample. In other words, obtain the integer $j_1$ for which $x_{j_1+n/2-1} - x_{j_1}$ reaches a minimum, with $1 \leq j_1 \leq n/2 + 1$.

2. Next, considering only the data in that interval, find the smallest interval that contains $n/4$ points, i.e., obtain the integer $j_2$ for which $x_{j_2+n/4-1} - x_{j_2}$ reaches a minimum, with $j_1 \leq j_2 \leq n/4 + j_1$.

3. Then continue this procedure to iteratively find minimum intervals of $n/8$ points, $n/16$ points, etc., with each interval lying within the preceding interval, until obtaining an interval with only two points, $x_{j_m}$ and $x_{j_m+1}$.

4. Finally, the mode is estimated by the mean of those values, $\widehat{M} = (x_{j_m} + x_{j_m+1})/2$.

Fig. 1 illustrates how to compute the HSM. The pseudocode of the Appendix calculates $\widehat{M}$ recursively. The HSM is the estimator of Robertson and Cryer [5] with the function $k(n)$ of their algorithm equal to the integer part of $(n+1)/2$. Sager [6] provided multivariate analogues to the estimators of Robertson and Cryer [5].

Fig. 1 here

The estimate $\widehat{M}$ tends to be close to the mode because minimizing the interval widths with a constant number of data points per interval is the same as maximizing the empirical density. To determine how close the HSM for an ideal sample is to the mode of its distribution $F$, let $x_i$ be the $[(i - \frac{1}{2})/n)]$th quantile of $F$ for each $i$ and compute $\widehat{M}$ for the quantiles instead of random numbers, using $n$ that is high enough that $x_1 \leq \widehat{M} \leq x_n$. In this case, since $F$ has only one mode, it follows that $x_{j_m} \leq M \leq x_{j_m+1}$ and hence $|\widehat{M} - M| \leq x_{j_m+1} - x_{j_m}$ and $\lim_{n \to \infty} |\widehat{M} - M| = 0$. For observed data, $|\widehat{M} - M|$ will not necessarily be this small, but $\widehat{M}$ will still be close to $M$ insofar as $x_i$ is close to the $[(i - \frac{1}{2})/n)]$th quantile for each $i$.

Like certain other estimators, $\widehat{M}$ uses counts within intervals to estimate the maximum density. However, the definition of $\widehat{M}$ in terms of varying interval widths according to a decreasing number of points per interval has advantages over estimation methods that are based on a constant interval width $w$, such as the midpoint of the modal interval and density estimation, in which case $w$ is a smoothing parameter. (These methods are described in detail in the next subsection.) First, $\widehat{M}$ does not require the arbitrary selection of an interval width, a problem posed by constant-interval methods. If the interval is too large, then the mode cannot be precisely located, leading to a high bias, but if the interval is too small, then it will be likely that the interval with the highest empirical frequency does not contain the mode, leading to a high variance in the estimates. $\widehat{M}$ avoids these issues by beginning with a large interval and progressively reducing the interval size to locate the mode more precisely. In addition, $\widehat{M}$ is less sensitive to outliers, as demonstrated in Section 3.

## 2.2 Other robust estimators of the mode

While the HSM is based on taking successive half samples, two other estimators of the mode are based on only taking the first half sample, the half of the points in the sample with the shortest range. The first of these, called the *shorth* [7], is the mean of those points and the second, called the one-dimensional *location of the least median of squares* (LMS) [8], is the midpoint of the highest and lowest of the points in the shortest half. More precisely, given $x_1 \leq x_2 \leq \ldots \leq x_n$ and the integer $m$ that minimizes $x_{m+h} - x_m$, the shorth is the mean of $x_m, x_{m+1}, \ldots, x_{m+h}$, and the location of the LMS is $(x_{m+h} - x_m)/2$, where $h = n/2$ for even $n$ or $h = (n-1)/2$ for odd $n$. Rousseeuw and Leroy [9] suggested the location of the LMS as an estimator of the mode; Bickel [10] also considered the shorth as a possible estimator of the mode.

There are other nonparametric mode estimators that are insensitive to outliers and to changes in the distribution. Perhaps the simplest is the midpoint of the estimated *modal interval*, the interval of some fixed width $w$ that includes the maximum number of data points of a sample [11, 1]. (Wegman [12] proved that this simple estimator of the mode has strong consistency.) The *half-range mode* (HRM) is based on finding modal intervals that are subsets of other modal intervals, where each $w$ is equal to half the difference between the largest and smallest values of the previous modal interval, starting with the entire sample. Thus, the HRM is similar to the HSM, except that it uses half-ranges instead of half-samples [10]. An older estimator of the mode is the maximum of an empirical density function [11, 13, 14], called the *empirical probability density function mode* (EPDFM). The one studied in Section 3 is that of [15], the value of $x$ that maximizes

$$\widehat{f}(x) = \frac{1}{nh\sqrt{2\pi}} \sum_{i=1}^{n} \exp\left[-\frac{1}{2}\left(\frac{x - x_i}{h}\right)^2\right], \tag{3}$$

where $\sigma$ is the lower of the standard deviation and the normal-consistent median absolute deviation from the median (MAD) and $h = 0.9\,\sigma/n^{1/5}$ is the smoothing parameter. This is a modification of an empirical density function recommended by Silverman [16], who used the interquartile range instead of the more robust MAD. The EPDFM can be seen as a generalization of the midpoint of the modal interval [11], with $h$ taking the place of $w$.

Each parametric mode (PM) estimator of [15] is the maximum of a probability density function with parameters estimated from the data. The PM studied in Section 3 is found by applying the power transform to the data that brings the data as close as possible to the normal distribution, estimating the mean and standard deviation of the transformed data with the sample median and MAD, and computing the mode of the original data under the assumption that the transformed data is normally distributed. A potential drawback is that this PM is not invariant to shifts in the sample: if $x_0$ is added to all values of the sample, the PM of the new sample will not be $x_0$ plus the PM of the old sample. If this is a problem, PM can be modified by including a shift parameter, as was recommended for cases in which the data can be negative [15]. Such a parameter is not added here since optimizing two transformation parameters (the shift parameter and the scaling exponent) requires much more computation time than only optimizing one transformation parameter.

## 3 Comparisons between mode estimators

### 3.1 Breakdown point

The robustness of an estimator of location can be quantified by its breakdown point, the lowest proportion of outliers needed to make the estimates useless. Mathematically, for a sample of $n$ values $\{x_1, x_2, \ldots, x_{n-\nu}, y_1, y_2, \ldots, y_\nu\}$ with each $x_i$ drawn from the distribution $F$, the *finite-sample breakdown point* $\varepsilon_n^*$ of an estimator $\widehat{\mu}$ is the minimum value of $\varepsilon_n \equiv \nu/n$ for which there exists some set of real numbers $(y_i)_{i=1}^{\nu}$ that will cause $\widehat{\mu}$ to be unbounded [1]. The *infinite-sample breakdown point*, defined by $\varepsilon^* \equiv \lim_{n \to \infty} \varepsilon_n^*$, will be used to compare the robustness of the mean, the median, and several estimators of the mode.

The sample mean, $\left(\sum_{i=1}^{n-\nu} x_i + \sum_{i=1}^{\nu} y_i\right)/n$, has a finite-sample breakdown point of $\varepsilon_n^* = 1/n$ since the single outlier $y_1$ can be made sufficiently large to bring the mean past any value. Thus, the infinite-sample breakdown point for the mean is $\varepsilon^* = 0$ and the mean is not robust. The median, on the other hand, is very robust, since at least $\nu = n/2$ sufficiently high values of $y_i$ are needed to make the median much higher than the maximum value of $x_i$. It follows that the finite-sample and infinite-sample breakdown points are 1/2. The trimmed means, as compromises between the mean and median, have breakdown points between 0 and 1/2, depending on how they are trimmed [1]. Quantiles other than the median also have intermediate breakdown points.

The robustness of the midpoint of the estimated modal interval depends on the distribution $F$ and the size $w$ of the interval. If $f$ is the probability density of $F$, $M$ is the only mode of $F$, and $\widetilde{M}(w)$ is the midpoint of the modal interval of $F$, then the number of values of $x_i$ in the modal interval is approximately $N(w) = (n - \nu)P(w)$, where

$$P(w) = \int_{\widetilde{M}(w)-w/2}^{\widetilde{M}(w)+w/2} f(x)\,dx$$

is the probability that $x_i$ lies in the interval. The estimated modal interval of the sample $\{x_1, x_2, \ldots, x_{n-\nu}, y_1, y_2, \ldots, y_\nu\}$ can be made much higher than all values of $x_i$ if $\nu \geq N(w)$. Hence, the breakdown points are $\varepsilon_n^* = N(w)/n$ and $\varepsilon^* = P(w)/[P(w) + 1]$. It is evident that these breakdown points increase as the interval width, $w$, increases, but higher values of $w$ imply that the mode, $M$, is estimated less closely since $|\widetilde{M}(w) - M| \leq w/2$. For example, an infinite-sample breakdown point as high as $\varepsilon^* = 1/2$ can be achieved for distributions that have zero density outside some interval, but only when $w$ is so large that $P(w) = 1$ and the modal interval includes all possible values of $x_i$, rendering the midpoint of the estimated modal interval a very poor estimator of the mode. Therefore, there is a trade-off between the robustness and precision of the estimated modal interval.

The estimator of the mode given by Eq. (2) is not robust since $M_{p,k}^*$ will greatly exceed all values of $x_i$ if $\nu \geq k + 1$ and if $y_1, y_2, \ldots, y_\nu$ are sufficiently high and sufficiently close. This can be seen from $\lim_{y_1 \to y_2 \to \cdots \to y_\nu \to \infty} M_{p,k}^* = \infty$. Therefore, $\varepsilon_n^* = (k+1)/n$ and $\varepsilon^* = 0$ for this estimator.

The HSM $\widehat{M}$, the estimator of the mode described in Section 2, will be bounded as long as less than half of the values of the sample are contamination values. To see this, let $\nu < n/2$ and $x_1 \leq x_2 \leq \ldots \leq x_{n-\nu}$. It follows that all intervals with $n/2$ points contain at least one value $x_i$ with $1 \leq i \leq n - \nu$. Let $[X_{\min}, X_{\max}]$ be the smallest closed interval that contains $n/2$ values of the sample. Thus, $\exists i \in \{1, 2, 3, \ldots, n - \nu\}$ such that $x_i \in [X_{\min}, X_{\max}]$. Since $n - \nu > n/2$ the smallest interval with $n/2$ points cannot be greater than the range of non-contamination values, i.e., $X_{\max} - X_{\min} \leq x_{n-\nu} - x_1$. Therefore, all values in that interval are bounded from below and from above:

$\forall\, x \in [X_{\min}, X_{\max}] : x_1 - (x_{n-\nu} - x_1) \leq x \leq x_{n-\nu} + (x_{n-\nu} - x_1)$. By the definition of the HSM, $\widehat{M} \in [X_{\min}, X_{\max}]$, implying that the estimate of the mode is likewise bounded: $x_1 - (x_{n-\nu} - x_1) \leq \widehat{M} \leq x_{n-\nu} + (x_{n-\nu} - x_1)$. The HSM can only be unbounded when the assumption that $\nu < n/2$ is relaxed, yielding $\nu > n/2$ and $\varepsilon_n \geq 1/2$. Hence, the finite-sample breakdown point of $\widehat{M}$ is $\varepsilon_n^* = 1/2$ and its infinite-sample breakdown point is $\varepsilon^* = 1/2$. In this respect, the HSM is as robust as the median. However, the HSM and median differ in the two measures of robustness of Section 3.2.

The shorth and the location of the LMS have the same breakdown point as the HSM since they are also based on the shortest half-sample. However, as estimators of the mode, they are highly biased for asymmetric distributions [10], which is confirmed in Section 3.3.

A less-biased estimator of the mode with $\varepsilon^* = 1/2$ is the HRM [10], and the PM of Section 2.2 also has a breakdown point of 1/2 [15]. The breakdown point of the EPDFM has not been determined analytically, but the simulations of Section 3.3 suggest that it is less than 1/2. (Additional simulations indicate that the infinite-sample breakdown point is close to 25% for the normal distribution and 40% for the lognormal distribution for the smoothing parameter computed as per Section 2.2 (R. Frühwirth, unpublished)). The HSM is much faster computationally than these three estimators, as seen in Section 3.4.

## 3.2 Sensitivity curves

Influence functions are often used to give the asymptotic effect of a single contamination point on an estimator for a particular distribution [7], but they are difficult or impossible to derive for the HSM and many of the other estimators of Section 2. Consequently, this section instead has *stylized sensitivity curves* (SSCs), which describe the effect of the addition of a single point on an estimator, given quantile approximations to the expected order statistics of a particular distribution. The SSC of estimator $T_n$ is defined as

$$S_n(x) = n\left[T_n(x_1, x_2, \ldots, x_{n-1}, x) - T_{n-1}(x_1, x_2, \ldots, x_{n-1})\right], \tag{4}$$

where $x$ is the location of a contamination point and the quantile

$$x_i = F^{-1}\left(\frac{i - 1/2}{n - 1}\right) \tag{5}$$

is defined in terms of a distribution $F$ [7]. The three distributions of Fig. 2 were used as $F$: a normal distribution, which has zero skewness and kurtosis, a lognormal distribution, which has positive skewness and kurtosis, and a Pareto distribution, which has infinite skewness and kurtosis. <span style="float:right">Fig. 2 here</span>

The SSCs of the estimators of Section 2 for $n = 100$ are displayed with that of the sample median in Figs. 3–10. (The sample median was chosen for comparison to estimators of the mode since, like the mode and unlike many other measures of location, the median is familiar to most researchers and is easily understood. The sample median is also the most widely used robust estimator. Bickel and Lehmann [2] similarly compared estimators of different location parameters.) Based on their influence function analogs [17], two measures of robustness have been defined in terms of an SSC: the finite-sample *rejection point* $\rho_n^*$, the smallest value of $|x|$ for which the addition of $x$ has no effect on the estimator, and the finite-sample *gross-error sensitivity* $\gamma_n^*$, the maximum possible effect of the addition of $x$ on the estimator [10]. That is,

$$\rho_n^* \equiv \inf\{r > 0 : S_n(x) = 0 \text{ if } |x| > r\} \qquad (6)$$

and

$$\gamma_n^* \equiv \sup_x |S_n(x)|. \qquad (7)$$

As seen in Fig. 3, the median has an infinite rejection point for all three distributions since the addition of any point shifts the median, but, except for $x$ very close to the uncontaminated median, it is shifted the same amount for all $x$ on either side of the median. The amount shifted is small, so the gross-error sensitivity is also small. The SSC of the standard PM (Fig. 4) reflects the fact that this estimator is not robust to outliers: both the rejection point and the gross-error sensitivity are infinite. Figs. 5–8 are similar in that their gross-error sensitivities are moderate and their rejection points are finite for the normal distribution, but infinite for at least one of the two skewed distributions. The HRM and HSM stand out from the other robust estimators, having high gross-error sensitivities and finite rejection points for all three distributions: the maxima of the SSCs in Figs. 9–10 are greater than those of any of the other figures, but the SSCs always go to zero for sufficiently small and sufficiently large values of $x$. In other words, the HRM and HSM completely reject all outliers, but are sensitive to contamination near the mode. As the sample size grows, their rejection points decrease and their gross-error sensitivities increase [10]. <span style="float:right">Figs. 3–10 here</span>

### 3.3 Simulations

The HSM is even more robust to outliers than the median when they are either much greater than or much less than all values of $x_i$ and when $\nu < n/2$. This is because $\widehat{M}$ is not biased by outliers of this type, but the median of a sample with these outliers is $x_{(n+1)/2}$, which is higher than $x_{(n-\nu+1)/2}$, the median of the sample without outliers, assuming $x_1 \leq x_2 \leq \ldots \leq x_{n-\nu} \leq y_1 \leq y_2 \leq \ldots \leq y_\nu$.

The robustness of the HSM relative to the median and the robust estimators of the mode of Section 2.2 was examined using simulated samples of $n$ random numbers each, with $n = 100$, 500, and 1000. For each sample, $(1-\varepsilon)n$ random numbers were drawn from a distribution $F$ and $\varepsilon n$ random numbers were drawn from an outlier distribution, where $\varepsilon$ is the portion of contamination. The contamination distribution associated with distribution $F$ is a normal distribution with a mean equal to the 99.99th percentile of $F$ and with a standard deviation equal to 1% of the interquartile range of the main distribution divided by the interquartile range of the standard normal distribution. 10,000 samples were generated from each distribution $F$ of Fig. 2 and at each value of $n$, with $\varepsilon = 0$, 0.1, 0.2, 0.3, and 0.4, and the median and mode estimates were computed for each sample. Tables 1–9 display the bias, standard error, and root mean square error (RMSE) of the estimates for the simulations. The standard error is the standard deviation of estimates. The bias is the mean error of the estimates, and the RMSE is the square root of the mean of the squared errors, where the error of a mode estimate is the estimate minus the mode of $F$ and the error of the sample median is the estimate minus the median of $F$. Thus, both the bias and standard error of an estimator contribute to its RMSE. (The bias and RMSE of the sample median would generally be higher for asymmetric distributions if its error were defined as the estimate minus the mode, treating the sample median as an estimator of the mode. We instead compare mode estimators to the sample median as an estimator of the median; such comparisons between estimators of different measures of location are meaningful since measures of location are translation invariant [2], according to Eq. (1).) The bias and RMSE of the estimates for $n$=500 are displayed graphically in Fig. 11.  <span style="float:right">Fig. 11 here</span>

Some trends in bias emerge from Tables 1–3. First, the bias in the median is much more sensitive to contamination than the bias in the mode estimators. For the normal distribution, the PM is the most resistant to high levels of contamination. For the lognormal and Pareto distributions, the bias in the HRM and HSM is approximately independent of the level of contamination, due to the finite rejection point of these estimators (Section 3.2), whereas the bias in the shorth and LMS increases as the contamination level increases. Even at zero contamination, the shorth and LMS are highly biased for the lognormal and Pareto distributions, as noted in Section 3.1. Overall, the mode estimators of lowest bias are the shorth, LMS, and EPDFM for the normal distribution with low levels of contamination, the PM for the lognormal distribution and for the Pareto distribution with small sample sizes, and the PM, HRM, and HSM for the Pareto distribution with large sample sizes.

[Tables 1–3 here]

Similar trends appear in the standard errors (Tables 4–6). The standard error of the PM tends to be less sensitive to the level of contamination than the other estimators. Among estimators of the mode, the PM, shorth, and EPDFM generally have the lowest standard errors for the normal and lognormal distributions. The low standard error of the EPDFM for the normal distribution at high levels of contamination and large sample sizes is deceptive since the mode that is consistently found is the mode from the contamination distribution, as is clear from Tables 1 and 7; see the Appendix on algorithms of the EPDFM. For the Pareto distribution, the PM has the lowest standard error given small sample sizes and the EPDFM, HRM, and HSM have the lowest standard errors given larger sample sizes.

[Tables 4–6 here]

Since the RMSE is a measure of performance that takes both bias and standard error into account, the findings of Tables 7–9 will be discussed with the conclusions.

[Tables 7–9 here]

## 3.4 Computation time

The amount of time it takes to compute an estimate becomes important when data must be simulated or resampled, as in bootstrapping and in permutation tests. The times depend on the particular algorithms and implementations of the estimators, but the relative times between estimators are similar for our *Mathematica* and MATLAB implementations, as seen in Figs. 12–13. The computation time increases with the sample size, but is comparatively independent of the distribution for all estimators except the HRM. For the sample sizes and distributions considered, the threshold of 1 second in the *Mathematica* implementations (Fig. 12) can be used to differentiate between fast estimators (GM, shorth, LMS, HSM) and slow estimators (PM, HRM, EPDFM). <span style="float:right">[Figs. 12–13 here]</span>

# 4 Applications

## 4.1 Robust initialization of the M-estimator

The M-estimator [7] is a popular robust location estimator which has not been considered so far in this paper. It can be computed via an iterated reweighted least-squares estimator, and the final result may depend on the initial location and scale. Rousseuw [9] has proposed to use the median and the normal-consistent median absolute deviation from the median (MAD) as initial estimates of location and scale, respectively. Another robust scale estimate is the length of the shorth (i.e. of the shortest half sample) [18]. It has also been considered as an initial scale estimate for the M-estimator [19]. It turns out that in some cases, such as the contaminated normal sample used above, better starting values can be obtained by using the HSM as the initial location and the normal-consistent half-width at half-maximum (HWHM) as the initial scale. The HWHM can be computed from the empirical probability density function (see Section 2.2) in the following way:

1. Compute the HSM $\widehat{M}$.

2. Compute the EPDF and look for the local maximum nearest to $\widehat{M}$.

3. Move out from this local maximum on each side until the EPDF drops to half of its value at the maximum.

4. Compute the half-width at the half-maximum and divide it by 1.1774. The result is a normal-consistent estimate of the scale.

Six different scenarios of initializing the M-estimator have been compared. They are summarized in Table 10. All initial scale estimates have been made normal-consistent by applying appropriate scale factors. <span style="float:right">[Table 10 here]</span>

Table 11 shows the RMSE of the M-estimator of the location, for normal samples with different levels of contamination. The samples have been generated in exactly the same way as in the simulation study reported above. The M-estimator uses Huber's

$\psi$-function, defined by

$$\psi(t) = \begin{cases} t, & \text{if } |t| \leq c \\ c\,\text{sgn}(t), & \text{if } |t| > c \end{cases}$$

The robustness constant $c$ has been set to 1.5.

The results show that at low contamination level the Median/MAD initialization is to be preferred, but that at high contamination level the HSM/HWHM is a better way of initializing the M-estimator. With the contaminated samples the length of the shorth always gives worse results than either the MAD or the HWHM.

Table 11 here

## 4.2 Finding the signal vertex in high-energy proton-proton collisions

The Large Hadron Collider (LHC) at CERN in Geneva is scheduled to go into operation in the year 2007 [20]. In the LHC, proton bunches will be accelerated to high energies (up to 7 Teraelectronvolts per beam). Once they have reached their final energy, they will be brought to collision at certain points along the accelerator. At the points of collision, huge detectors will record the particles produced by the collisions.

A crossing of two bunches is called an *event*. Because of the high density of the proton bunches, up to twenty individual collisions take place simultaneously in each event. The event rate is very high, about 40 Megahertz. As the maximum recording rate is only a hundred Hertz, the rate is reduced by a multi-level trigger which selects the events which are the most interesting ones from the point of view of physics. In a selected event there is a *signal vertex* which is the collision point of the interesting process, plus up to twenty background vertices which are randomly distributed over the entire length of the crossing region (about 30 centimeters).

In the course of the data analysis it is important to quickly find the signal vertex, using the tracks of the particles produced in the event. The average transverse momentum of the particles emerging from the signal vertex is considerably higher than the one of the background vertices. This is illustrated by Figure 14, showing a histogram of the position $z$ of all tracks exceeding some momentum threshold, weighted by the transverse momentum. The distribution clearly has many modes, each one corresponding to a collision. The task of finding the signal vertex is now equivalent to finding the highest mode of the distribution of $z$, taking into account the weight of each observation.

To this end, the HSM was generalized in two ways. First, each observations is counted by its weight. (In this data set the weight is the transverse momentum, but other weights may apply to other problems. In general, the weight of an observation has the same effect as a proportional number of observations of the same value, so that an observation with twice the weight of another counts as if twice as many observations were made at its value. The use of weights is a simple way to map two variables to a single variable, so that a one-dimensional estimator of the mode can be applied.) The content of an interval is therefore the sum of the weights of the observations in the interval. Second, each successive interval is required to contain a certain fraction $p$ of the content of its predecessor, not necessarily equal to one half. The resulting estimator is called the *fraction-of-sample mode with weights* (FSMW). If all weights are equal to 1, it is the estimator of Robertson and Cryer [5] with the function $k(n) = p \cdot n$. The computation of the FSMW is only slightly more complicated than the standard HSM and can be implemented in not more than 20 lines of MATLAB code.

Alternative mode estimators in this context are the mode of the weighted frequency histogram (HISTMW) and the mode of the weighted empirical density function (EPDFMW). All three estimators considered have one tunable parameter: the fraction $p$ in the case of the FSMW, the bin size in the case of the HISTMW, and the smoothing parameter $h$ in the case of the EPDFMW (see Eq. (3)).

The data sample used in this study consists of 981 simulated events from the CMS experiment [21]. As the true signal vertex is known in each event, the parameters of the three methods could be optimized. The optimization was done by searching for the parameter values giving the smallest mean-squared deviation of the estimated signal vertex from the true signal vertex. This procedure yielded an optimal bin size of 0.008 cm for the HISTMW estimator, an optimal smoothing parameter $h = 0.018$ cm for the EPDFMW, and an optimal fraction $p = 0.06$ for the FSMW.

Table 12 shows the bias, the standard deviation and the RMSE of the three mode estimators with respect to the true signal vertex. The precision of the three estimators is comparable, the FSMW having the smallest RMSE, followed by the HISTMW and the EPDFMW.

### 4.3 City sizes

Our final example of data analysis with the mode is the application of the HSM to the 1930 sizes of 49 cities. As seen from Fig. 15, the distribution is highly asymmetric; this leads to very different values of the mean, median, and mode. Following Davison and Hinkley [22], it is assumed that the cities compose an IID sample so that inferences can be made using bootstrapping methods (this assumption is not needed to estimate the mode). The standard error of the HSM was estimated as the standard deviation of estimates for 2000 samples of 49 cities each, drawn from the original data with replacement. The same bootstrap samples were used to estimate the quartiles of HSM estimates, using the bias-correction method of Efron [23]. The same procedure was applied to the sample mean and sample median, for comparison. Table 13 indicates that the standard error and interquartile range of the HSM are much lower than those of the mean and median. Using the fact that the HSM estimate is equal to three city sizes and greater than two city sizes, the *modal skewness* [10] estimate is

$$\widehat{\alpha} = 1 - 2\left\{\left(\#\left[x = \widehat{M}\right]/2 + \#\left[x < \widehat{M}\right]\right)/n\right\} = 1 - 2\{(3/2 + 2)/49\} \approx 86\%, \quad (8)$$

reflecting the asymmetry of this data. That estimate is close to the modal skewness of the lognormal distribution of Fig. 2: $\alpha = 1 - 2F(M) \approx 89\%$.   Table 11 here    Fig. 14 here

## 5 Discussion and Conclusions

Whether the HSM or the sample median is a better estimator of location depends on the application. When the standard error or confidence intervals of the estimates are needed but only one sample is available, an advantage of the median is that the distribution of its estimates can be computed if the distribution of the data is known. However, the distribution of mode estimates can be easily computed from a single sample using resampling (Section 4) and bootstrap estimates can be used to make inferences related to the mode. For example, the hypothesis that the modes of two populations differ can be tested by replacing the sample mean with the HSM in bootstrap methods of detecting a difference in means [24, 22]. To test for differences in the modes of marginal distributions of many variables, the HSM can be used with multiple comparison procedures [25]. If a sample is known to be from a distribution with low skewness and kurtosis, and if it is known that there are few outliers, then the median might be more appropriate. In other cases, the HSM may be a better estimator of location. It is also useful as a highly robust starting point for iterative location estimators such as the M-estimator.

The RMSE, as recorded in Tables 7-9, can be used to evaluate the performance of the robust estimators of the mode under different conditions. For the normal distribution, the estimators with the lowest RMSE are the PM for no contamination, the shorth or EPDFM for moderate contamination, and the LMS for high contamination. For the lognormal and Pareto distributions, the PM consistently performs best, except in the Pareto case with a high level of contamination and a large sample size, when the HRM and HSM perform better. Since the HSM is orders of magnitude faster than the HRM (Figs. 12-13), while having very similar statistical properties, the HSM should be used instead of the HRM in most situations.

Considering only the robust mode estimators that are also fast, as defined in Section 3.4, the shorth and LMS perform better than the HSM for the normal distribution, but the HSM performs best in the Pareto case for all sample sizes and levels of contamination. For the lognormal distribution, the shorth and LMS have lower RMSE for low contamination and $n = 20$, but the HSM has the lowest RMSE of the fast estimators when the sample size is higher or when the contamination is moderate to high. Due to its finite rejection point, the HSM always has a lower bias than the shorth and LMS for the lognormal and Pareto distributions (Section 3.1 and Tables 2–3), which matters for applications in which a low bias is more important than a low variance.

The finite rejection point of the HSM also enables its use in testing for differences in the mode. Although the mode itself does not satisfy condition (iv) of Section 1, the infinite rejection points of most other estimators of the mode for lognormal or Pareto data can lead to the problem of two samples with mode estimates that appear to be statistically significant due to differences in the tails rather than in the central parts of the distributions. The HSM, on the other hand, is not influenced by tails and is thus better-suited for inferences about modes for data that strongly deviates from normality. (Many of the mode estimators have finite rejection points for normal data and thus may be suitable if the distribution is only slightly skewed.) In biomedicine and other disciplines, inferences about values of high probability density are often of more interest than inferences about parameters that are influenced by values of low probability density [25].

In conclusion, any of the following conditions suggests the use of the HSM:

1. The distribution is extremely asymmetric with a high level of contamination with outliers.

2. A fast estimator is needed and the distribution is at least moderately asymmetric.

3. An estimator with a finite rejection point is needed, as in cases of testing for differences between modes or in cases in which the bias must be independent of the level of contamination, and the distribution is at least moderately asymmetric.

When none of these conditions is satisfied, the PM or another estimator of the mode might have better performance. However, the HSM performs relatively well under a wide range of conditions (Fig. 11; Tables 1–9), so it can serve as a general-purpose mode estimator when the distribution or level of contamination is unknown.

## 6 Acknowledgements

We thank Susanna Cucciarelli and Wolfgang Waltenberger for providing us with the data used in the signal vertex study. We also thank the anonymous reviewers for comments that helped us improve this paper.

## Appendix: Algorithms and Implementations

The estimators of the mode described in Section 2 were implemented in both *Mathematica* and MATLAB. Equivalent algorithms were used in *Mathematica* and MATLAB, except for the two estimators requiring function optimization over a real argument: the PM and the EPDFM. The *Mathematica* implementation of PM followed the algorithm of Bickel [15] for finding $\beta_0$, the transformation exponent that maximizes $R(\beta)$, the normality of the transformed data, quantified as a robust measure of correlation between the transformed data and normal expected order statistics. The MATLAB implementation of PM followed a maximization algorithm that is faster computationally, apparently without significant loss of precision:

1. Let $\beta_L = -2.9$, $\beta_R = 4.1$, and $\Delta\beta = 0.15$.

2. Let $\beta_0 = \arg\max R(\beta)$ for $\beta \in \{\beta_L, \beta_L + \Delta\beta, \ldots, \beta_R - \Delta\beta, \beta_R\}$.

3. If $\beta_L + 6\,\Delta\beta < \beta_0 < \beta_R - 6\,\Delta\beta$, then go to Step 4.

    If $\beta_0 \leq \beta_L + 6\,\Delta\beta$, set $\beta_L \longleftarrow \beta_L - 6\,\Delta\beta$ and go to Step 2.

    If $\beta_0 \geq \beta_R - 6\,\Delta\beta$, set $\beta_R \longleftarrow \beta_R + 6\,\Delta\beta$ and go to Step 2.

4. Let $\beta_1 = \beta_0 - \Delta\beta$, $\beta_2 = \beta_0$, and $\beta_3 = \beta_0 + \Delta\beta$.

5. If $|R(\beta_3) - R(\beta_1)| \leq 10^{-4}$, then go to Step 6. Otherwise, set $\Delta\beta \longleftarrow \Delta\beta/2$, set $\beta_0 = \arg\max R(\beta)$ for $\beta \in \{\beta_2 - \Delta\beta, \beta_2, \beta_2 + \Delta\beta\}$, and go to Step 4.

6. The PM is the argument of the vertex of the parabola that intersects each point in $\{R(\beta_1), R(\beta_2), R(\beta_3)\}$.

The *Mathematica* implementation of the EPDFM, following that used in [15], often settles on a local maximum of the density for the mode estimate rather than the global maximum when the portion of extreme contamination is large. The MATLAB implementation consistently uses the global maximum of the density, estimating the mode by $\arg\max_x \widehat{f}(x)$. In the simulations reported in [15], the failure of the *Mathematica* implementation to find the global maximum unexpectedly worked to the advantage of the EPDFM since the algorithm was more robust to outliers than the more theoretically transparent MATLAB implementation used in Section 3.3, above. The poor performance of the MATLAB version of the EPDFM at high levels of contamination results from high density estimates at the contamination regions.

The following algorithm, called **HalfSampleMode**, estimates the mode from a sample of $n$ real numbers, ordered such that $x_1 \leq x_2 \leq \ldots \leq x_n$. The theory behind **HalfSampleMode** is explained in Section 2. *Mathematica* and MATLAB versions of this algorithm were used to compute $\widehat{M}$ in Sections 3–4, with the MATLAB version used to generate Tables 1–9 and Figs. 3–11. The MATLAB version of the algorithm is iterative and is about 2/3 faster than a recursive MATLAB version, with little sample size dependence. Both an iterative version in S code and the recursive version in *Mathematica* code can be found at www.davidbickel.com. The S version has been tested in both S-PLUS and R. R is similar to the S statistical computing environment, but R runs on more operating systems than S-PLUS does, and R is freely available from www.r-project.org.

**HalfSampleMode**$(n, (x_1, x_2, \ldots, x_n))$

1. If there is only one number in the sample ($n = 1$), then that number is the estimate of the mode. Return $x_1$ and stop.

2. If there are two numbers in the sample ($n = 2$), then the estimate of the mode is the mean of those numbers. Return $(x_1 + x_2)/2$ and stop.

3. If there are three numbers in the sample ($n = 3$), then the estimate of the mode is the mean of the two numbers that are closest together: if $x_2 - x_1 < x_3 - x_2$, then return $(x_1 + x_2)/2$ and stop, but if $x_2 - x_1 < x_3 - x_2$, then return $(x_2 + x_3)/2$ and stop. Otherwise, if $n=3$ and $x_2 - x_1 = x_3 - x_2$, then the mode is estimated to be $x_2$; return $x_2$ and stop.

4. If there are four or more numbers in the sample ($n \geq 4$), then the algorithm will be recursively applied to about half of the sample. Set the minimum interval width to the range of the sample: $w_{\min} = x_n - x_1$. This will be replaced by the width of the smallest interval that contains approximately half of the values of the sample. Set $N$ equal to the smallest integer greater than or equal to $n/2$. (This algorithm could be generalized by setting $N$ equal to the smallest integer greater than or equal to $\alpha n$, where $0 < \alpha < 1$, but estimators have not been studied with $\alpha \neq 1/2$.) For $i = 1, 2, 3, \ldots, n - N + 1$, do the following two steps:

   (a) Set the width $w = x_{i+N-1} - x_i$.

   (b) If $w < w_{\min}$, then set $w_{\min} = w$ and set $j = i$.

   Return **HalfSampleMode**$(N, (x_j, x_{j+1}, \ldots, x_{j+N-1}))$. Note that this recursive call reduces the problem of estimating the mode of the sample of size $n$ to that of estimating the mode of the reduced sample of size $N$ that has the smallest range of values.

# Table captions

Table 1:   Bias of estimators for normal distribution

Table 2:   Bias of estimators for lognormal distribution

Table 3:   Bias of estimators for Pareto distribution

Table 4:   Standard error of estimators for normal distribution

Table 5:   Standard error of estimators for lognormal distribution

Table 6:   Standard error of estimators for Pareto distribution

Table 7:   RMSE of estimators for normal distribution

Table 8:   RMSE of estimators for lognormal distribution

Table 9:   RMSE of estimators for Pareto distribution

Table 10: Six scenarios for the initial location and scale of the M-estimator

Table 11: RMSE of the M-estimator for normal distribution with six different initial values

Table 12: Bias, standard deviation and RMSE of the estimated mode with respect to the true signal vertex.

Table 13: Estimated mean, median, and mode of the city size data of Fig. 14 (in thousands of people)

| $n$ | $\varepsilon$ | Median | PM | Shorth | LMS | EPDFM | HRM | HSM |
|---|---|---|---|---|---|---|---|---|
| 20 | 0% | 0.004 | −0.041 | 0.006 | 0.006 | 0.002 | −0.002 | 0.000 |
|  | 10% | 0.135 | −0.080 | −0.011 | −0.011 | −0.005 | −0.014 | −0.009 |
|  | 20% | 0.314 | −0.135 | 0.001 | 0.000 | 0.003 | −0.002 | −0.001 |
|  | 30% | 0.562 | −0.320 | 0.050 | 0.044 | 0.015 | 0.081 | 0.071 |
|  | 40% | 0.956 | −0.198 | 0.788 | 0.618 | 1.053 | 1.139 | 1.083 |
| 100 | 0% | 0.001 | −0.010 | −0.001 | 0.000 | −0.003 | −0.003 | −0.001 |
|  | 10% | 0.138 | −0.084 | −0.002 | −0.003 | −0.001 | 0.003 | 0.002 |
|  | 20% | 0.315 | −0.196 | −0.002 | −0.002 | −0.005 | −0.004 | 0.000 |
|  | 30% | 0.563 | −0.429 | −0.002 | −0.002 | −0.001 | 0.012 | −0.002 |
|  | 40% | 0.961 | −0.311 | 0.895 | 0.659 | 3.645 | 1.972 | 1.119 |
| 500 | 0% | −0.001 | −0.005 | −0.001 | −0.001 | −0.003 | −0.003 | −0.002 |
|  | 10% | 0.139 | −0.071 | 0.000 | 0.000 | −0.001 | 0.005 | 0.003 |
|  | 20% | 0.319 | −0.210 | −0.001 | −0.002 | −0.001 | −0.002 | 0.002 |
|  | 30% | 0.565 | −0.459 | −0.001 | −0.001 | 1.759 | 0.010 | −0.001 |
|  | 40% | 0.965 | −0.324 | 0.975 | 0.708 | 3.717 | 3.500 | 1.160 |
| 1000 | 0% | −0.001 | −0.003 | 0.000 | 0.000 | 0.001 | 0.002 | 0.001 |
|  | 10% | 0.140 | −0.065 | 0.001 | 0.002 | −0.002 | −0.004 | −0.005 |
|  | 20% | 0.318 | −0.214 | 0.000 | 0.000 | −0.001 | −0.003 | −0.005 |
|  | 30% | 0.566 | −0.463 | 0.000 | 0.001 | 3.718 | 0.014 | 0.002 |
|  | 40% | 0.967 | −0.323 | 0.986 | 0.713 | 3.718 | 3.701 | 1.172 |

Table 1: Bias of estimators for normal distribution

| $n$ | $\varepsilon$ | Median | PM | Shorth | LMS | EPDFM | HRM | HSM |
|---|---|---|---|---|---|---|---|---|
| 20 | 0% | 0.123 | 0.370 | 0.851 | 0.913 | 1.000 | 0.706 | 0.724 |
|  | 10% | 0.536 | 0.284 | 0.949 | 1.045 | 1.148 | 0.749 | 0.757 |
|  | 20% | 1.204 | 0.202 | 1.122 | 1.286 | 1.379 | 0.815 | 0.826 |
|  | 30% | 2.384 | −0.011 | 1.395 | 1.712 | 1.705 | 0.870 | 0.883 |
|  | 40% | 5.148 | 0.020 | 1.921 | 2.648 | 2.236 | 0.964 | 0.981 |
| 100 | 0% | 0.026 | 0.044 | 0.612 | 0.691 | 0.635 | 0.335 | 0.342 |
|  | 10% | 0.429 | 0.025 | 0.731 | 0.856 | 0.766 | 0.348 | 0.350 |
|  | 20% | 1.045 | 0.015 | 0.907 | 1.118 | 0.955 | 0.371 | 0.379 |
|  | 30% | 2.114 | −0.031 | 1.181 | 1.589 | 1.230 | 0.390 | 0.391 |
|  | 40% | 4.526 | −0.061 | 1.695 | 2.688 | 1.678 | 0.429 | 0.428 |
| 500 | 0% | 0.002 | 0.007 | 0.549 | 0.629 | 0.421 | 0.167 | 0.170 |
|  | 10% | 0.412 | −0.005 | 0.673 | 0.800 | 0.533 | 0.174 | 0.180 |
|  | 20% | 1.029 | −0.016 | 0.854 | 1.070 | 0.699 | 0.185 | 0.187 |
|  | 30% | 2.078 | 0.101 | 1.135 | 1.552 | 0.945 | 0.196 | 0.197 |
|  | 40% | 4.444 | −0.079 | 1.647 | 2.679 | 8.892 | 0.240 | 0.212 |
| 1000 | 0% | 0.001 | 0.003 | 0.539 | 0.619 | 0.348 | 0.127 | 0.125 |
|  | 10% | 0.411 | −0.011 | 0.665 | 0.791 | 0.454 | 0.136 | 0.133 |
|  | 20% | 1.023 | −0.023 | 0.846 | 1.062 | 0.607 | 0.134 | 0.136 |
|  | 30% | 2.075 | 0.158 | 1.128 | 1.545 | 0.840 | 0.143 | 0.144 |
|  | 40% | 4.447 | −0.081 | 1.643 | 2.680 | 103.688 | 0.179 | 0.159 |

Table 2: Bias of estimators for lognormal distribution

| $n$ | $\varepsilon$ | Median | PM | Shorth | LMS | EPDFM | HRM | HSM |
|---|---|---|---|---|---|---|---|---|
| 20 | 0% | 0.793 | 0.246 | 1.296 | 1.726 | 1.622 | 0.689 | 0.686 |
| | 10% | 2.384 | 0.202 | 1.632 | 2.329 | 2.031 | 0.761 | 0.778 |
| | 20% | 6.064 | 0.113 | 2.252 | 3.632 | 2.789 | 0.832 | 0.843 |
| | 30% | 18.645 | $-0.075$ | 3.721 | 7.195 | 4.533 | 0.912 | 0.922 |
| | 40% | 276.965 | $-0.260$ | 10.743 | 30.988 | 14.057 | 1.043 | 1.039 |
| 100 | 0% | 0.136 | 0.124 | 1.055 | 1.541 | 1.078 | 0.292 | 0.296 |
| | 10% | 1.247 | 0.074 | 1.314 | 2.071 | 1.316 | 0.300 | 0.304 |
| | 20% | 3.501 | 0.004 | 1.760 | 3.133 | 1.709 | 0.320 | 0.325 |
| | 30% | 9.391 | $-0.199$ | 2.689 | 5.858 | 2.488 | 0.347 | 0.344 |
| | 40% | 40.350 | $-0.377$ | 5.686 | 19.193 | 4.874 | 0.371 | 0.373 |
| 500 | 0% | 0.026 | 0.110 | 1.011 | 1.508 | 0.849 | 0.145 | 0.146 |
| | 10% | 1.101 | 0.039 | 1.263 | 2.041 | 1.041 | 0.151 | 0.153 |
| | 20% | 3.185 | $-0.026$ | 1.685 | 3.071 | 1.353 | 0.158 | 0.161 |
| | 30% | 8.468 | $-0.243$ | 2.540 | 5.675 | 1.954 | 0.168 | 0.169 |
| | 40% | 33.368 | $-0.324$ | 5.124 | 17.808 | 3.698 | 0.180 | 0.182 |
| 1000 | 0% | 0.012 | 0.108 | 1.005 | 1.503 | 0.775 | 0.111 | 0.113 |
| | 10% | 1.083 | 0.032 | 1.258 | 2.037 | 0.954 | 0.116 | 0.117 |
| | 20% | 3.144 | $-0.037$ | 1.674 | 3.062 | 1.238 | 0.121 | 0.122 |
| | 30% | 8.347 | $-0.244$ | 2.518 | 5.643 | 1.793 | 0.126 | 0.129 |
| | 40% | 32.624 | $-0.272$ | 5.054 | 17.629 | 3.398 | 0.134 | 0.137 |

Table 3: Bias of estimators for Pareto distribution

| $n$ | $\varepsilon$ | Median | PM | Shorth | LMS | EPDFM | HRM | HSM |
|---|---|---|---|---|---|---|---|---|
| 20 | 0% | 0.271 | 0.423 | 0.474 | 0.474 | 0.439 | 0.563 | 0.565 |
| | 10% | 0.284 | 0.464 | 0.460 | 0.462 | 0.400 | 0.575 | 0.579 |
| | 20% | 0.305 | 0.490 | 0.441 | 0.445 | 0.372 | 0.581 | 0.590 |
| | 30% | 0.336 | 0.481 | 0.561 | 0.537 | 0.380 | 0.817 | 0.787 |
| | 40% | 0.398 | 0.562 | 1.476 | 1.214 | 1.634 | 1.801 | 1.773 |
| 100 | 0% | 0.126 | 0.256 | 0.264 | 0.271 | 0.295 | 0.395 | 0.390 |
| | 10% | 0.133 | 0.288 | 0.252 | 0.263 | 0.265 | 0.400 | 0.397 |
| | 20% | 0.140 | 0.274 | 0.234 | 0.253 | 0.230 | 0.409 | 0.402 |
| | 30% | 0.156 | 0.224 | 0.216 | 0.247 | 0.213 | 0.472 | 0.413 |
| | 40% | 0.189 | 0.194 | 1.497 | 1.134 | 0.481 | 1.895 | 1.754 |
| 500 | 0% | 0.056 | 0.128 | 0.147 | 0.156 | 0.200 | 0.277 | 0.274 |
| | 10% | 0.059 | 0.163 | 0.138 | 0.151 | 0.183 | 0.287 | 0.283 |
| | 20% | 0.064 | 0.151 | 0.129 | 0.147 | 0.156 | 0.291 | 0.291 |
| | 30% | 0.071 | 0.110 | 0.116 | 0.144 | 1.861 | 0.346 | 0.298 |
| | 40% | 0.086 | 0.084 | 1.512 | 1.109 | 0.001 | 0.878 | 1.748 |
| 1000 | 0% | 0.040 | 0.093 | 0.115 | 0.124 | 0.171 | 0.242 | 0.240 |
| | 10% | 0.042 | 0.128 | 0.109 | 0.121 | 0.155 | 0.247 | 0.246 |
| | 20% | 0.045 | 0.117 | 0.100 | 0.118 | 0.131 | 0.253 | 0.251 |
| | 30% | 0.050 | 0.084 | 0.088 | 0.113 | 0.001 | 0.295 | 0.260 |
| | 40% | 0.061 | 0.059 | 1.512 | 1.104 | 0.001 | 0.257 | 1.740 |

Table 4: Standard error of estimators for normal distribution

| $n$ | $\varepsilon$ | Median | PM | Shorth | LMS | EPDFM | HRM | HSM |
|---|---|---|---|---|---|---|---|---|
| 20 | 0% | 0.783 | 0.931 | 0.745 | 0.732 | 0.677 | 0.969 | 0.972 |
|  | 10% | 0.941 | 0.862 | 0.726 | 0.724 | 0.686 | 1.019 | 1.022 |
|  | 20% | 1.238 | 0.757 | 0.744 | 0.778 | 0.750 | 1.087 | 1.084 |
|  | 30% | 1.801 | 0.594 | 0.800 | 0.931 | 0.856 | 1.149 | 1.161 |
|  | 40% | 3.414 | 0.662 | 0.986 | 1.448 | 1.074 | 1.682 | 1.313 |
| 100 | 0% | 0.347 | 0.342 | 0.274 | 0.284 | 0.257 | 0.541 | 0.540 |
|  | 10% | 0.418 | 0.278 | 0.274 | 0.295 | 0.254 | 0.558 | 0.553 |
|  | 20% | 0.530 | 0.268 | 0.281 | 0.321 | 0.275 | 0.591 | 0.591 |
|  | 30% | 0.763 | 0.302 | 0.299 | 0.401 | 0.306 | 0.617 | 0.617 |
|  | 40% | 1.404 | 0.263 | 0.385 | 0.689 | 0.384 | 0.646 | 0.644 |
| 500 | 0% | 0.152 | 0.125 | 0.121 | 0.137 | 0.115 | 0.351 | 0.349 |
|  | 10% | 0.186 | 0.121 | 0.119 | 0.140 | 0.108 | 0.360 | 0.356 |
|  | 20% | 0.239 | 0.139 | 0.121 | 0.151 | 0.111 | 0.369 | 0.367 |
|  | 30% | 0.343 | 0.211 | 0.133 | 0.187 | 0.124 | 0.379 | 0.376 |
|  | 40% | 0.616 | 0.115 | 0.166 | 0.311 | 27.751 | 1.958 | 0.392 |
| 1000 | 0% | 0.108 | 0.085 | 0.088 | 0.102 | 0.086 | 0.296 | 0.292 |
|  | 10% | 0.132 | 0.090 | 0.087 | 0.105 | 0.080 | 0.302 | 0.296 |
|  | 20% | 0.169 | 0.108 | 0.086 | 0.112 | 0.079 | 0.303 | 0.304 |
|  | 30% | 0.240 | 0.164 | 0.093 | 0.136 | 0.084 | 0.317 | 0.313 |
|  | 40% | 0.435 | 0.083 | 0.117 | 0.222 | 27.481 | 1.603 | 0.328 |

Table 5: Standard error of estimators for lognormal distribution

| $n$ | $\varepsilon$ | Median | PM | Shorth | LMS | EPDFM | HRM | HSM |
|---|---|---|---|---|---|---|---|---|
| 20 | 0% | 2.632 | 0.748 | 0.802 | 1.169 | 1.052 | 0.838 | 0.806 |
|  | 10% | 4.762 | 0.696 | 1.097 | 1.886 | 1.443 | 0.959 | 0.970 |
|  | 20% | 9.931 | 0.680 | 1.741 | 3.832 | 2.399 | 1.085 | 1.143 |
|  | 30% | 45.704 | 0.652 | 3.897 | 11.647 | 5.404 | 1.199 | 1.262 |
|  | 40% | 4477.236 | 0.746 | 35.755 | 155.577 | 58.901 | 1.565 | 1.628 |
| 100 | 0% | 0.875 | 0.258 | 0.252 | 0.427 | 0.264 | 0.255 | 0.257 |
|  | 10% | 1.299 | 0.262 | 0.332 | 0.629 | 0.335 | 0.265 | 0.268 |
|  | 20% | 2.339 | 0.269 | 0.486 | 1.119 | 0.469 | 0.288 | 0.287 |
|  | 30% | 5.643 | 0.319 | 0.873 | 2.632 | 0.782 | 0.318 | 0.310 |
|  | 40% | 32.606 | 0.397 | 2.592 | 13.810 | 2.076 | 0.355 | 0.346 |
| 500 | 0% | 0.364 | 0.136 | 0.105 | 0.181 | 0.089 | 0.113 | 0.111 |
|  | 10% | 0.543 | 0.141 | 0.140 | 0.270 | 0.114 | 0.120 | 0.119 |
|  | 20% | 0.939 | 0.134 | 0.202 | 0.465 | 0.158 | 0.125 | 0.124 |
|  | 30% | 2.154 | 0.184 | 0.356 | 1.064 | 0.257 | 0.134 | 0.131 |
|  | 40% | 10.004 | 0.249 | 0.937 | 4.879 | 0.619 | 0.146 | 0.144 |
| 1000 | 0% | 0.257 | 0.101 | 0.074 | 0.128 | 0.058 | 0.084 | 0.083 |
|  | 10% | 0.383 | 0.107 | 0.098 | 0.191 | 0.074 | 0.090 | 0.088 |
|  | 20% | 0.644 | 0.098 | 0.139 | 0.321 | 0.100 | 0.092 | 0.091 |
|  | 30% | 1.499 | 0.148 | 0.246 | 0.745 | 0.165 | 0.097 | 0.097 |
|  | 40% | 6.793 | 0.221 | 0.647 | 3.359 | 0.400 | 0.105 | 0.104 |

Table 6: Standard error of estimators for Pareto distribution

| $n$ | $\varepsilon$ | Median | PM | Shorth | LMS | EPDFM | HRM | HSM |
|---|---|---|---|---|---|---|---|---|
| 20 | 0% | 0.271 | 0.425 | 0.474 | 0.474 | 0.439 | 0.563 | 0.565 |
|  | 10% | 0.314 | 0.471 | 0.460 | 0.462 | 0.400 | 0.575 | 0.579 |
|  | 20% | 0.437 | 0.509 | 0.441 | 0.445 | 0.372 | 0.581 | 0.590 |
|  | 30% | 0.655 | 0.578 | 0.563 | 0.539 | 0.380 | 0.821 | 0.791 |
|  | 40% | 1.035 | 0.596 | 1.673 | 1.363 | 1.944 | 2.131 | 2.078 |
| 100 | 0% | 0.126 | 0.256 | 0.264 | 0.271 | 0.295 | 0.395 | 0.390 |
|  | 10% | 0.191 | 0.300 | 0.252 | 0.263 | 0.265 | 0.400 | 0.397 |
|  | 20% | 0.345 | 0.337 | 0.234 | 0.253 | 0.230 | 0.409 | 0.402 |
|  | 30% | 0.584 | 0.483 | 0.216 | 0.247 | 0.213 | 0.472 | 0.413 |
|  | 40% | 0.980 | 0.366 | 1.744 | 1.311 | 3.677 | 2.735 | 2.081 |
| 500 | 0% | 0.056 | 0.128 | 0.147 | 0.156 | 0.200 | 0.277 | 0.274 |
|  | 10% | 0.152 | 0.178 | 0.138 | 0.151 | 0.183 | 0.287 | 0.283 |
|  | 20% | 0.325 | 0.258 | 0.129 | 0.147 | 0.156 | 0.291 | 0.291 |
|  | 30% | 0.570 | 0.472 | 0.116 | 0.144 | 2.561 | 0.346 | 0.298 |
|  | 40% | 0.969 | 0.334 | 1.799 | 1.316 | 3.717 | 3.609 | 2.098 |
| 1000 | 0% | 0.040 | 0.093 | 0.115 | 0.124 | 0.171 | 0.242 | 0.240 |
|  | 10% | 0.146 | 0.143 | 0.109 | 0.121 | 0.155 | 0.247 | 0.246 |
|  | 20% | 0.322 | 0.244 | 0.100 | 0.118 | 0.131 | 0.253 | 0.251 |
|  | 30% | 0.568 | 0.471 | 0.088 | 0.113 | 3.718 | 0.295 | 0.260 |
|  | 40% | 0.969 | 0.329 | 1.805 | 1.314 | 3.718 | 3.710 | 2.098 |

Table 7: RMSE of estimators for normal distribution

| $n$ | $\varepsilon$ | Median | PM | Shorth | LMS | EPDFM | HRM | HSM |
|---|---|---|---|---|---|---|---|---|
| 20 | 0% | 0.792 | 1.002 | 1.131 | 1.171 | 1.207 | 1.199 | 1.212 |
|  | 10% | 1.083 | 0.908 | 1.194 | 1.271 | 1.338 | 1.264 | 1.272 |
|  | 20% | 1.726 | 0.784 | 1.346 | 1.503 | 1.569 | 1.358 | 1.363 |
|  | 30% | 2.988 | 0.595 | 1.608 | 1.949 | 1.908 | 1.442 | 1.459 |
|  | 40% | 6.178 | 0.662 | 2.159 | 3.018 | 2.481 | 1.939 | 1.639 |
| 100 | 0% | 0.348 | 0.345 | 0.671 | 0.747 | 0.685 | 0.636 | 0.639 |
|  | 10% | 0.600 | 0.279 | 0.781 | 0.905 | 0.807 | 0.657 | 0.655 |
|  | 20% | 1.172 | 0.268 | 0.949 | 1.163 | 0.994 | 0.697 | 0.702 |
|  | 30% | 2.247 | 0.303 | 1.219 | 1.639 | 1.268 | 0.730 | 0.730 |
|  | 40% | 4.739 | 0.270 | 1.738 | 2.775 | 1.721 | 0.775 | 0.774 |
| 500 | 0% | 0.152 | 0.125 | 0.563 | 0.644 | 0.436 | 0.389 | 0.388 |
|  | 10% | 0.452 | 0.121 | 0.684 | 0.812 | 0.544 | 0.400 | 0.399 |
|  | 20% | 1.056 | 0.139 | 0.863 | 1.081 | 0.708 | 0.413 | 0.412 |
|  | 30% | 2.106 | 0.234 | 1.143 | 1.563 | 0.953 | 0.426 | 0.425 |
|  | 40% | 4.486 | 0.140 | 1.656 | 2.697 | 29.140 | 1.973 | 0.445 |
| 1000 | 0% | 0.108 | 0.085 | 0.546 | 0.628 | 0.359 | 0.322 | 0.317 |
|  | 10% | 0.431 | 0.090 | 0.670 | 0.798 | 0.461 | 0.331 | 0.325 |
|  | 20% | 1.037 | 0.110 | 0.850 | 1.068 | 0.612 | 0.331 | 0.333 |
|  | 30% | 2.089 | 0.228 | 1.132 | 1.551 | 0.845 | 0.347 | 0.345 |
|  | 40% | 4.468 | 0.116 | 1.647 | 2.689 | 107.268 | 1.613 | 0.364 |

Table 8: RMSE of estimators for lognormal distribution

| $n$ | $\varepsilon$ | Median | PM | Shorth | LMS | EPDFM | HRM | HSM |
|---|---|---|---|---|---|---|---|---|
| 20 | 0% | 2.749 | 0.788 | 1.524 | 2.085 | 1.933 | 1.085 | 1.058 |
|  | 10% | 5.325 | 0.725 | 1.966 | 2.996 | 2.491 | 1.224 | 1.243 |
|  | 20% | 11.636 | 0.690 | 2.846 | 5.280 | 3.679 | 1.367 | 1.420 |
|  | 30% | 49.361 | 0.657 | 5.388 | 13.690 | 7.054 | 1.506 | 1.563 |
|  | 40% | 4485.795 | 0.790 | 37.335 | 158.633 | 60.556 | 1.881 | 1.931 |
| 100 | 0% | 0.886 | 0.286 | 1.085 | 1.599 | 1.110 | 0.387 | 0.392 |
|  | 10% | 1.801 | 0.272 | 1.356 | 2.165 | 1.358 | 0.400 | 0.405 |
|  | 20% | 4.210 | 0.269 | 1.825 | 3.327 | 1.772 | 0.430 | 0.433 |
|  | 30% | 10.956 | 0.376 | 2.827 | 6.422 | 2.608 | 0.471 | 0.463 |
|  | 40% | 51.878 | 0.547 | 6.249 | 23.645 | 5.298 | 0.513 | 0.509 |
| 500 | 0% | 0.365 | 0.175 | 1.016 | 1.519 | 0.853 | 0.184 | 0.183 |
|  | 10% | 1.228 | 0.146 | 1.271 | 2.059 | 1.047 | 0.193 | 0.193 |
|  | 20% | 3.321 | 0.136 | 1.698 | 3.106 | 1.362 | 0.201 | 0.203 |
|  | 30% | 8.738 | 0.305 | 2.565 | 5.773 | 1.971 | 0.215 | 0.213 |
|  | 40% | 34.836 | 0.409 | 5.209 | 18.465 | 3.749 | 0.232 | 0.232 |
| 1000 | 0% | 0.257 | 0.148 | 1.008 | 1.509 | 0.777 | 0.139 | 0.140 |
|  | 10% | 1.149 | 0.111 | 1.262 | 2.045 | 0.957 | 0.146 | 0.147 |
|  | 20% | 3.209 | 0.105 | 1.680 | 3.078 | 1.242 | 0.152 | 0.152 |
|  | 30% | 8.480 | 0.285 | 2.530 | 5.692 | 1.800 | 0.159 | 0.162 |
|  | 40% | 33.323 | 0.351 | 5.096 | 17.947 | 3.422 | 0.170 | 0.172 |

Table 9: RMSE of estimators for Pareto distribution

| Method | Initial location | Initial scale |
|---|---|---|
| (a) | Sample mean | Sample standard deviation |
| (b) | Median | Median absolute deviation from the median |
| (c) | Median | Length of the shorth |
| (d) | Half-sample mode | Length of the shorth |
| (e) | Median | Half-width at half-maximum |
| (f) | Half-sample mode | Half-width at half-maximum |

Table 10: Six scenarios for the initial location and scale of the M-estimator

| $n$ | $\varepsilon$ | (a) | (b) | (c) | (d) | (e) | (f) |
|---|---|---|---|---|---|---|---|
| 20 | 0% | 0.226 | 0.236 | 0.225 | 0.225 | 0.231 | 0.231 |
| | 10% | 0.340 | 0.296 | 0.337 | 0.337 | 0.312 | 0.312 |
| | 20% | 0.704 | 0.489 | 0.614 | 0.614 | 0.533 | 0.533 |
| | 30% | 1.140 | 0.909 | 1.041 | 1.041 | 0.909 | 0.909 |
| | 40% | 1.512 | 1.492 | 1.496 | 1.496 | 1.433 | 1.433 |
| 100 | 0% | 0.103 | 0.106 | 0.105 | 0.105 | 0.103 | 0.103 |
| | 10% | 0.271 | 0.202 | 0.238 | 0.238 | 0.224 | 0.224 |
| | 20% | 0.668 | 0.435 | 0.540 | 0.540 | 0.471 | 0.471 |
| | 30% | 1.129 | 0.890 | 1.066 | 1.066 | 0.853 | 0.853 |
| | 40% | 1.507 | 1.512 | 1.506 | 1.506 | 1.478 | 1.438 |
| 500 | 0% | 0.045 | 0.047 | 0.045 | 0.045 | 0.045 | 0.045 |
| | 10% | 0.256 | 0.178 | 0.216 | 0.216 | 0.201 | 0.201 |
| | 20% | 0.663 | 0.427 | 0.524 | 0.524 | 0.454 | 0.454 |
| | 30% | 1.128 | 0.885 | 1.067 | 1.067 | 0.816 | 0.816 |
| | 40% | 1.506 | 1.510 | 1.503 | 1.503 | 1.453 | 1.335 |
| 1000 | 0% | 0.032 | 0.033 | 0.032 | 0.032 | 0.032 | 0.032 |
| | 10% | 0.254 | 0.175 | 0.209 | 0.209 | 0.198 | 0.198 |
| | 20% | 0.662 | 0.426 | 0.524 | 0.524 | 0.450 | 0.451 |
| | 30% | 1.128 | 0.884 | 1.064 | 1.064 | 0.807 | 0.809 |
| | 40% | 1.507 | 1.510 | 1.504 | 1.504 | 1.420 | 1.311 |

Table 11: RMSE of the M-estimator for normal distribution with six different initial values

| Estimator | Bias | St.Dev. | RMSE |
|---|---|---|---|
| FSMW | 0.098 | 1.538 | 1.542 |
| HISTMW | 0.114 | 1.617 | 1.621 |
| EPDFMW | 0.083 | 1.659 | 1.662 |

Table 12: Bias, standard deviation and RMSE of the estimated mode with respect to the true signal vertex.

| Estimator | Estimate | Std. Error | 1st Quartile | Median | 3rd Quartile |
|---|---|---|---|---|---|
| Sample Mean | 128 | 17 | 117 | 129 | 141 |
| Sample Median | 79 | 10 | 75 | 79 | 85 |
| HSM | 50 | 6 | 48 | 50 | 50 |

Table 13: Estimated mean, median, and mode of the city size data of Fig. 14 (in thousands of people).

# Figure captions

Fig. 1: Data points of each shortest half sample used to compute the HSM of 16 values generated from the lognormal distribution with mean and variance parameters of 1. The 16 values of the first iteration in (a) represents the original sample. The sample with the 8 of those values with the smallest range is displayed for the second iteration in (a) and on a smaller scale in (b). The third and fourth iterations have 4 and 2 points, respectively, and the points can be best distinguished in (b). The mean of the 2 points of the last iteration is 1.25, so that is the value of the HSM.

Fig. 2: Probability densities of the three distributions from which samples were drawn; each distribution has a mode of 1. The normal and lognormal distributions each have a mean parameter of 1 and a variance parameter of 1. (The normal distribution used in the sensitivity plots and simulations actually has a mean parameter of 6 in order to ensure the positive numbers needed for the PM, but the shift of the mean does not affect the bias or variance of the other estimators.) The Pareto distribution has a cutoff parameter of 1 and a shape parameter of 1/2; its PDF is $1/(2\,x^{3/2})$ for $x \geq 1$ and 0 for $x < 1$. Pareto distributions are sometimes used to model income. The median of the normal distribution is 1, the median of the lognormal distribution is 2.72, and the median of the Pareto distribution is 4.

Fig. 3: Sensitivity curve of the median for (a) the normal distribution, (b) the lognormal distribution, and (c) the Pareto distribution.

Fig. 4: Sensitivity curve of the standard PM for (a) the normal distribution, (b) the lognormal distribution, and (c) the Pareto distribution. The standard PM is the same as the PM described in Section 2.2, except that the sample mean and standard deviation are used instead of the median and median deviation. It has low variance in the absence of outliers (Bickel, 2001). However, as seen from the unbounded sensitivity curve, the standard PM, unlike the PM of Fig. 5, has a breakdown point of 0%.

Fig. 5:  Sensitivity curve of the PM described in Section 2.2 for (a) the normal distribution, (b) the lognormal distribution, and (c) the Pareto distribution. This estimator has a high breakdown point, so it was used in the simulations of Section 3.3.

Fig. 6:  Sensitivity curve of the shorth for (a) the normal distribution, (b) the lognormal distribution, and (c) the Pareto distribution.

Fig. 7:  Sensitivity curve of the LMS for (a) the normal distribution, (b) the lognormal distribution, and (c) the Pareto distribution.

Fig. 8:  Sensitivity curve of the EPDFM for (a) the normal distribution, (b) the lognormal distribution, and (c) the Pareto distribution.

Fig. 9:  Sensitivity curve of the HRM for (a) the normal distribution, (b) the lognormal distribution, and (c) the Pareto distribution.

Fig. 10: Sensitivity curve of the HSM for (a) the normal distribution, (b) the lognormal distribution, and (c) the Pareto distribution. Note the similarity to Fig. 9.

Fig. 11: Bias and RMSE of location estimators for samples of $n=500$ values from the normal, lognormal, and Pareto distributions, contaminated by outliers from an extreme normal distribution. The HRM is excluded since its bias and RMSE are close to those of the HSM. The RMSE scale is logarithmic, so the vertical axis begins at 0.001 in the normal case and at 0.01 in the lognormal and Pareto cases. Numeric values are recorded in Tables 1–9.

Fig. 12: Time to compute the Grenander mode $M^*_{2,3}$ (GM) and the six robust mode estimators of Section 2, using the *Mathematica* implementations: *Mathematica*, Version 4.1 for Macintosh (Wolfram Research, Inc.), running on a 533 MHz G4 PowerPC processor (Apple Computer, Inc.); see the Appendix for more information. The computation time is the mean number of seconds it took to compute an estimate from a sample of size $n$, out of 1000 random samples from the distributions of Fig. 2, with $n = 250, 500, 1000$. (For the HRM, only 100 samples were used, due to its slow computation time.)

Fig. 13: Time to compute the Grenander mode $M^*_{2,3}$ (GM) and the six robust mode estimators of Section 2, using the MATLAB implementations: MATLAB, Version

5.1 for Windows (The MathWorks, Inc.), running on a 600 MHz Athlon processor (Advanced Micro Devices, Inc.); see the Appendix for more information. The computation time is the mean number of seconds it took to compute an estimate from a sample of size $n$, out of 1000 random samples from the distributions of Fig. 2, with $n = 250, 500, 1000$. (For the HRM, PM, and EPDFM, only 200 samples were used, due to their slow computation times.)

Fig. 14: Frequency distribution of the track position $z$ for a typical event. The frequencies are weighted by the transverse track momenta. The small triangle indicates the position of the signal vertex. The data are simulated data from the CMS experiment.

Fig. 15: Empirical cumulative distribution function (CDF) of city sizes (data from Davison and Hinkley, 1997). The mode of a distribution is sometimes defined as the vertex of its CDF, the point below which the CDF is convex and above which it is concave (Dharmadhikari and Joag-dev, 1988).

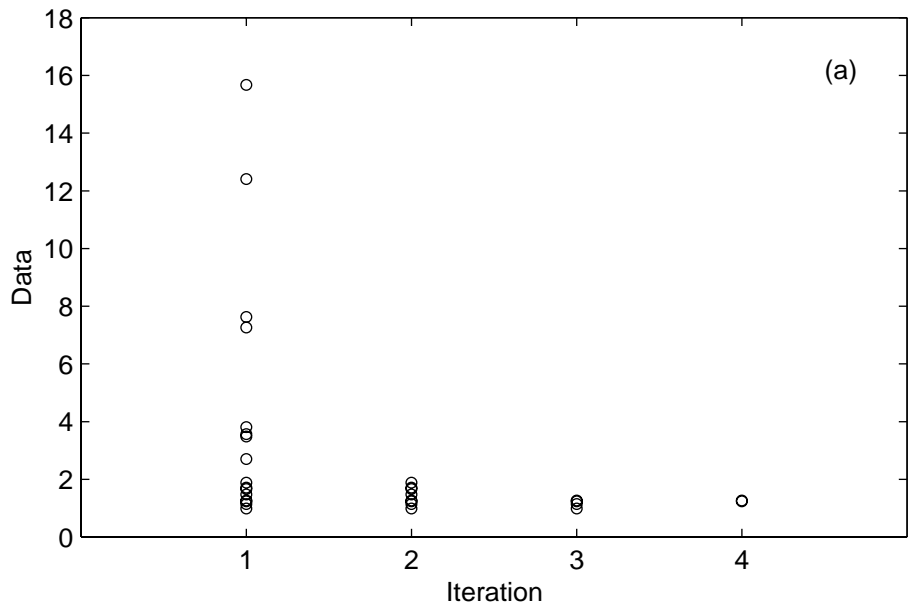

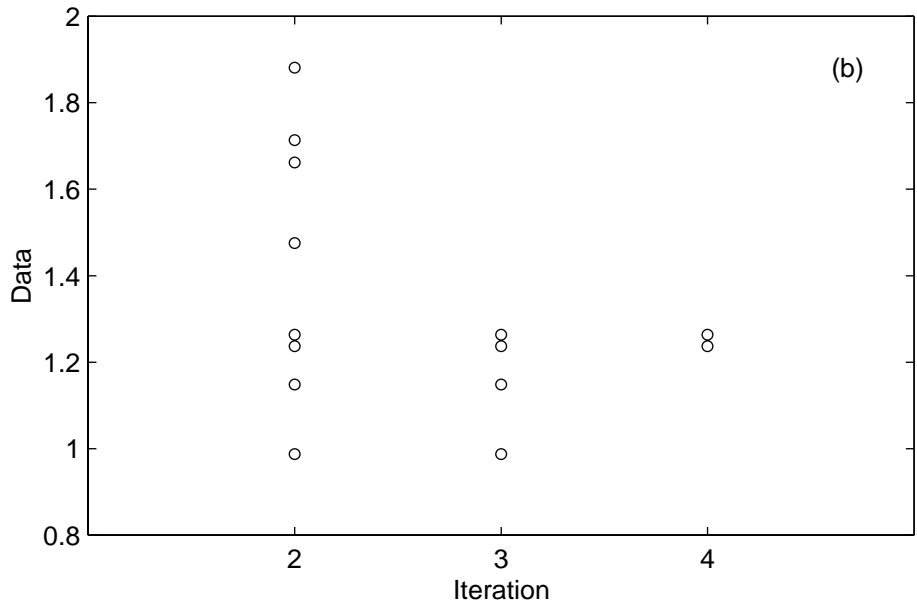

Figure 1

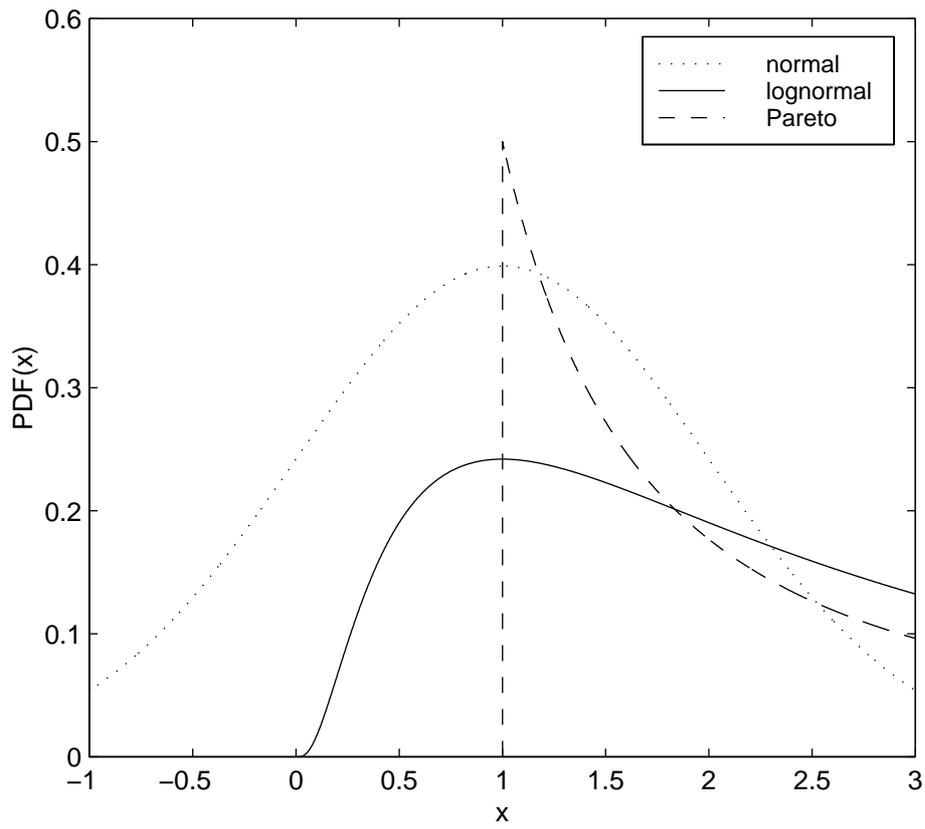

Figure 2

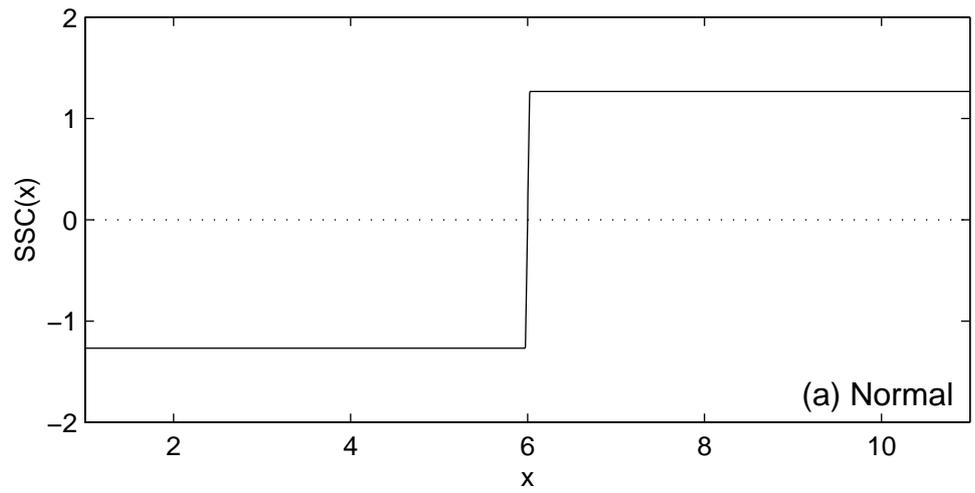

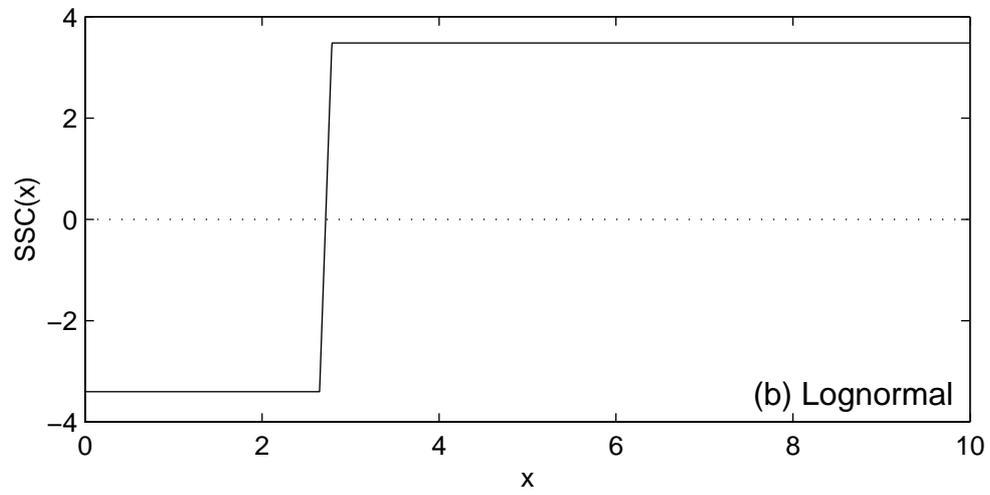

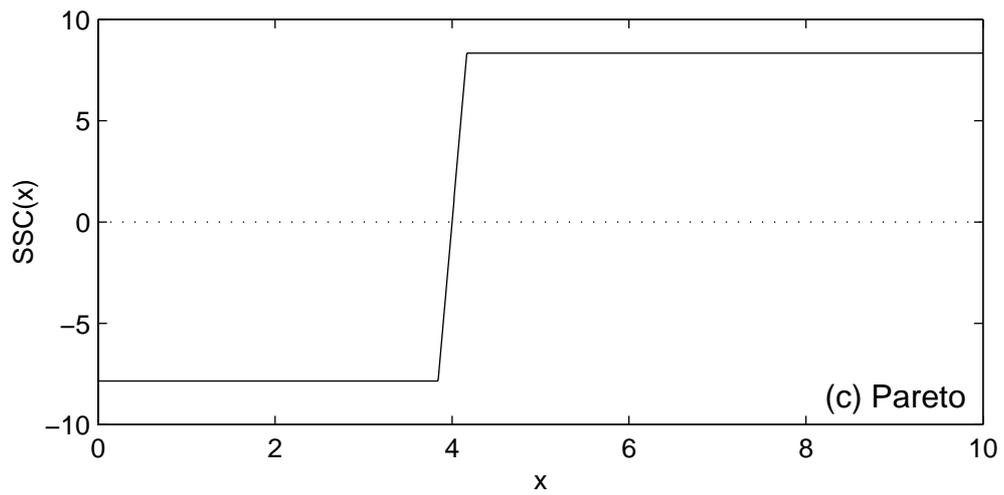

Figure 3

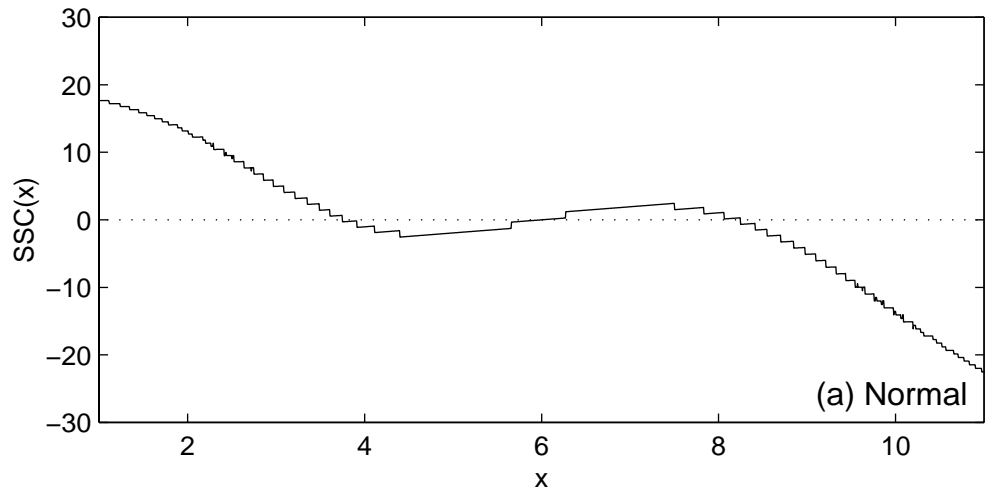

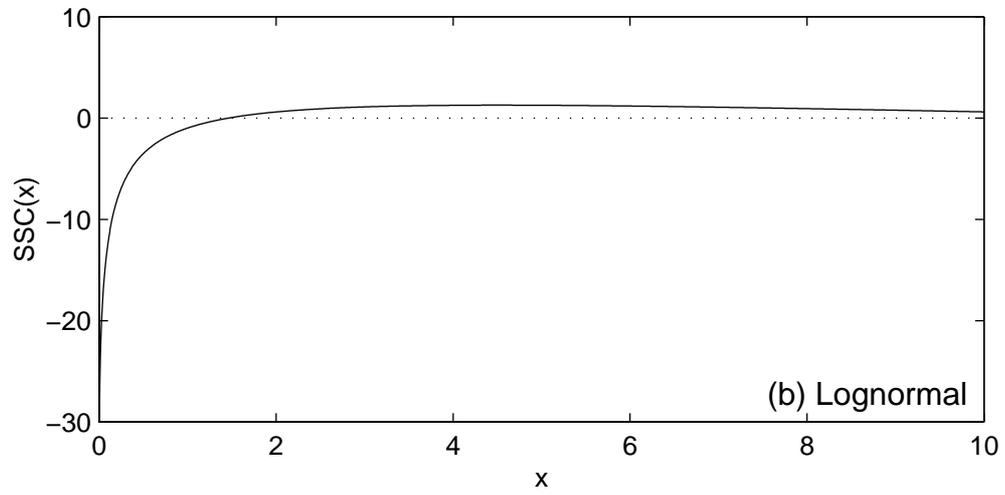

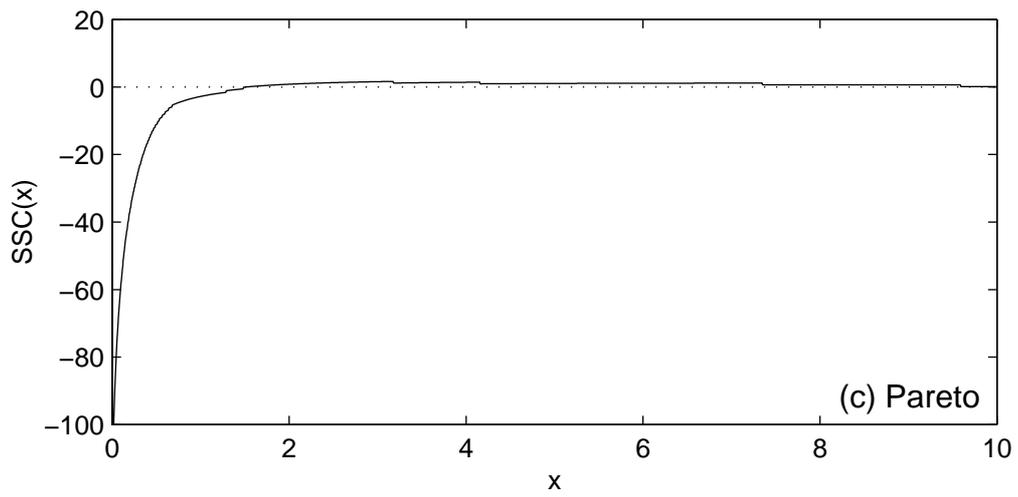

Figure 4

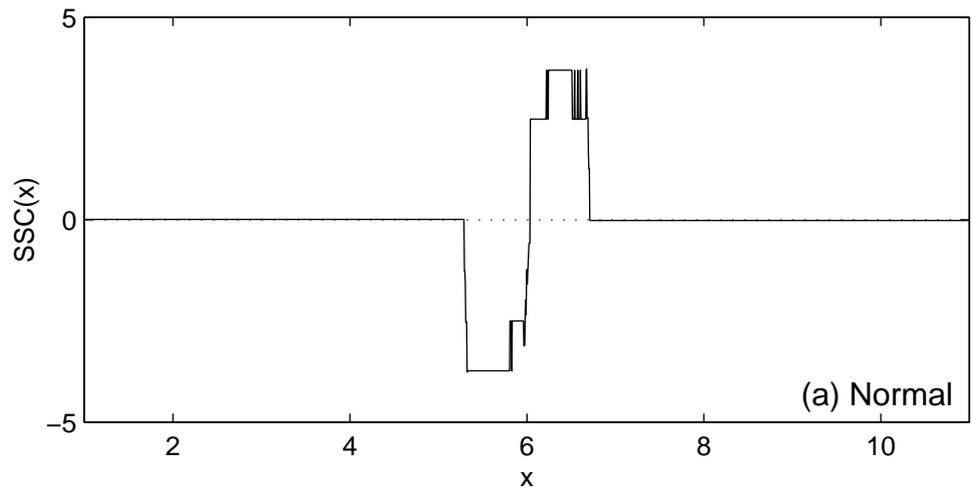
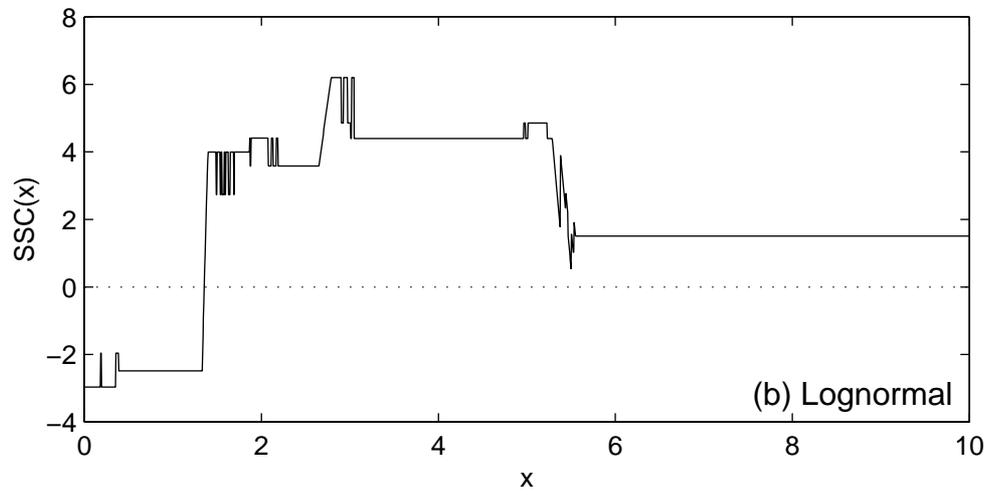
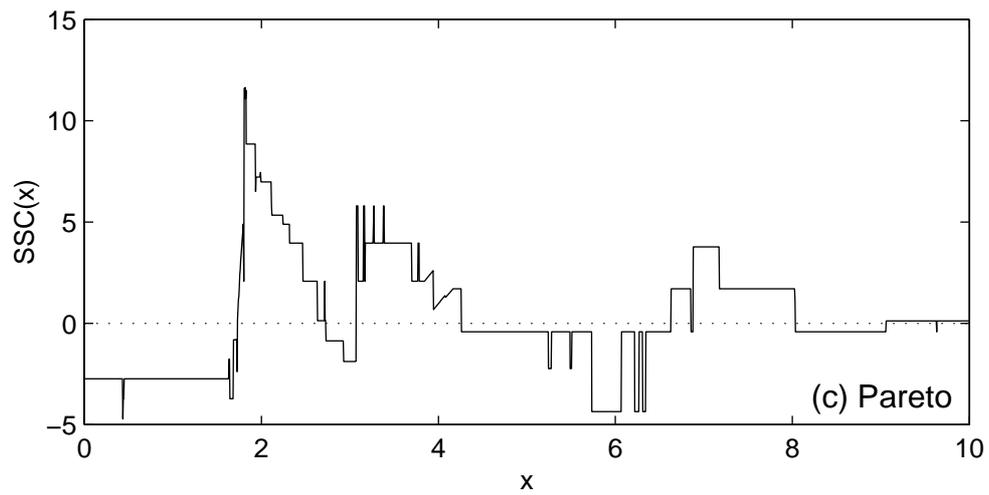

Figure 5

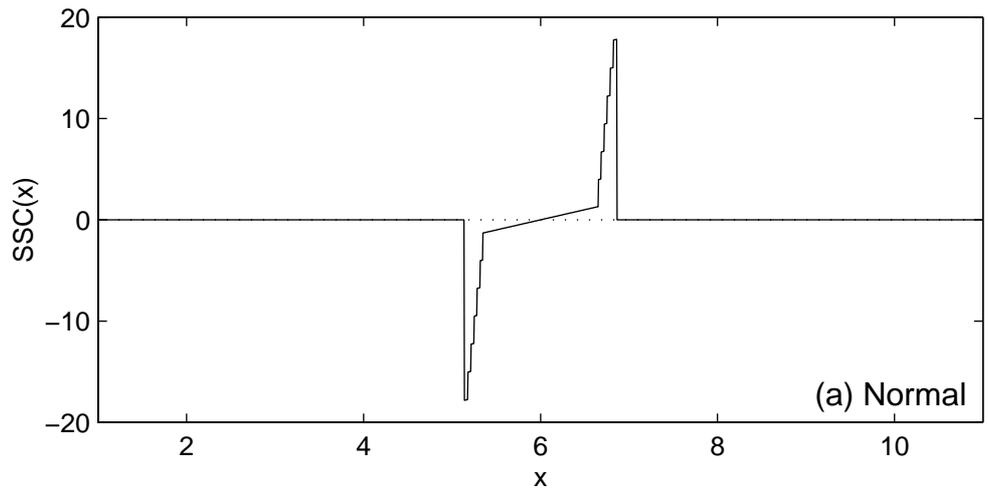
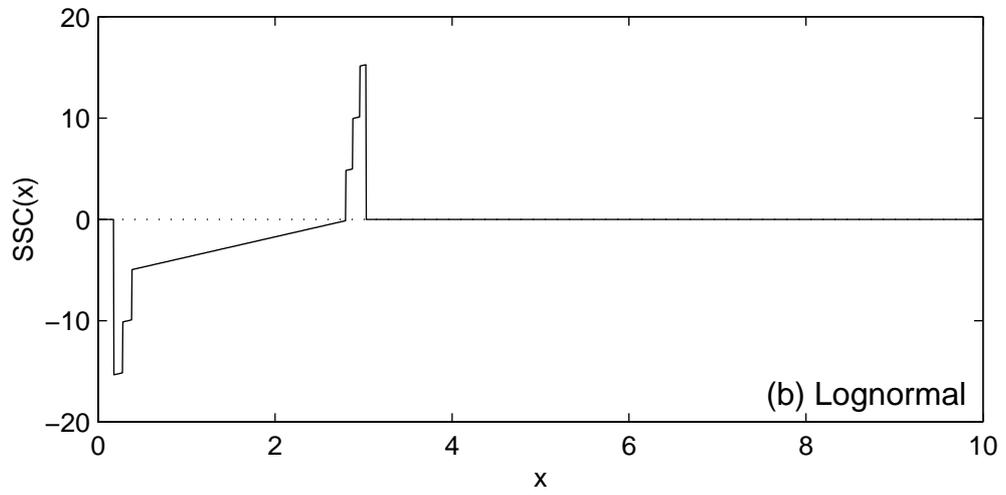
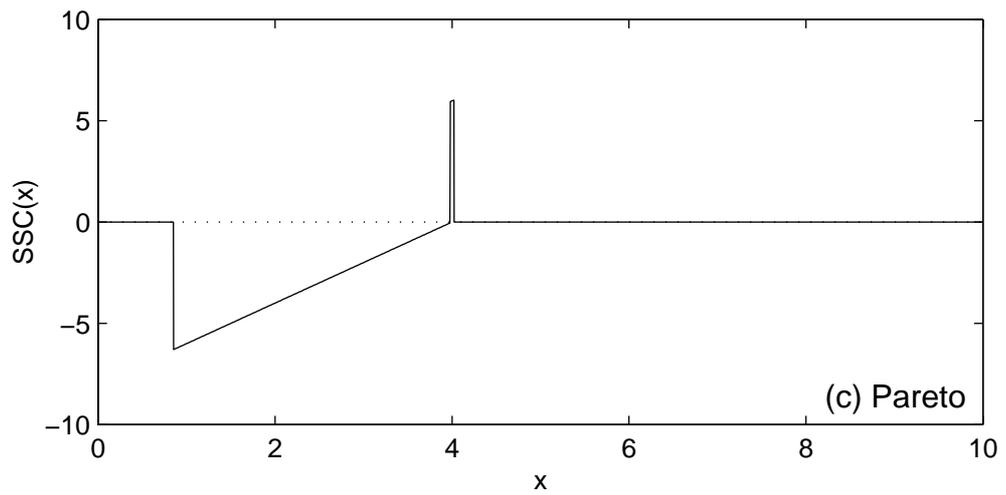

Figure 6

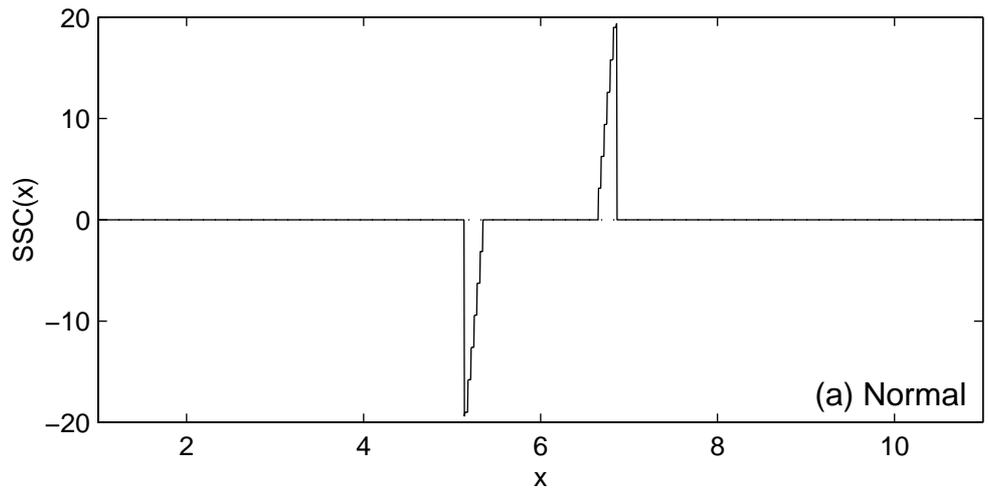
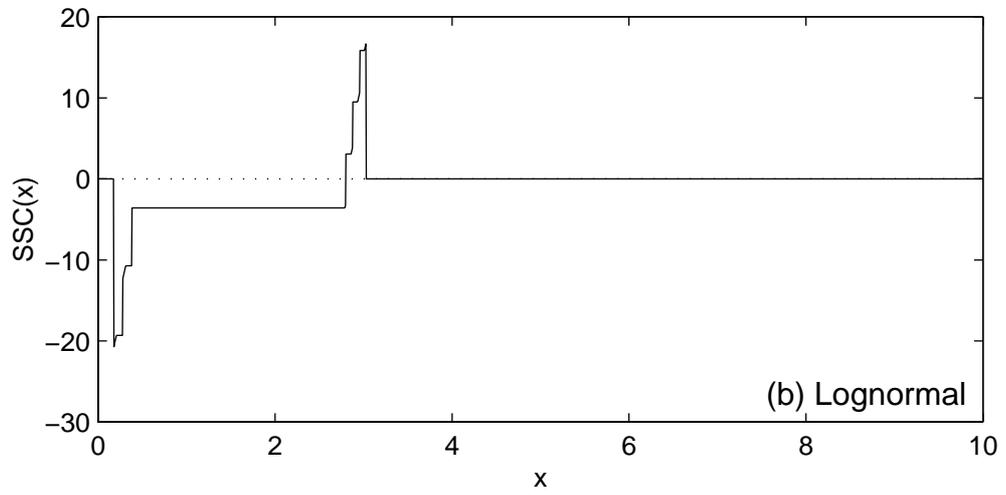
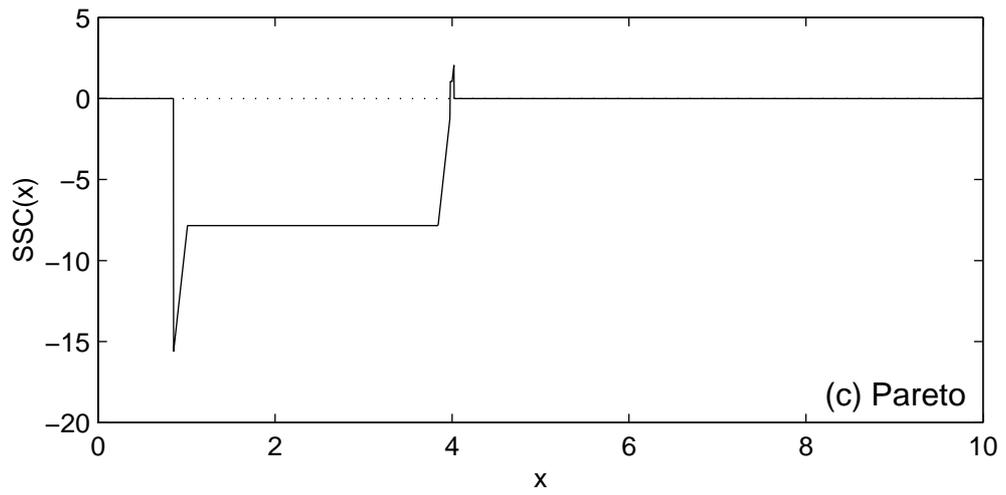

Figure 7

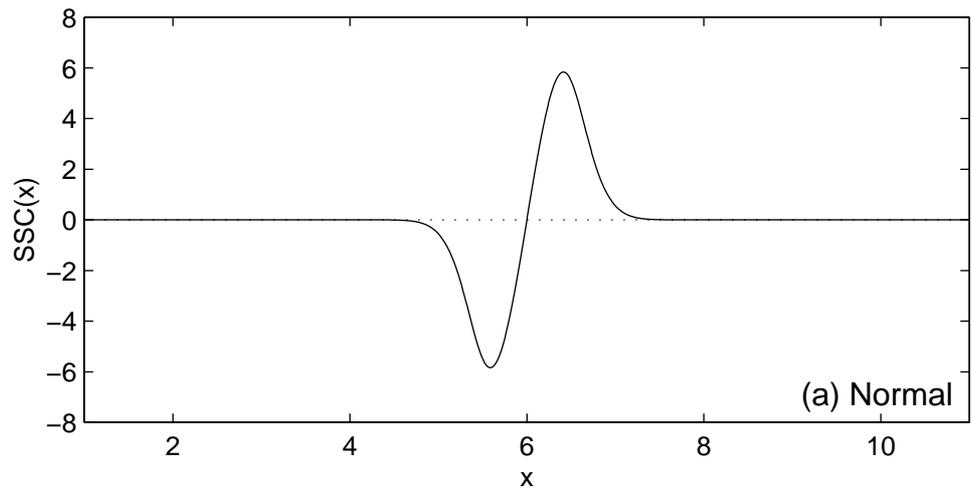

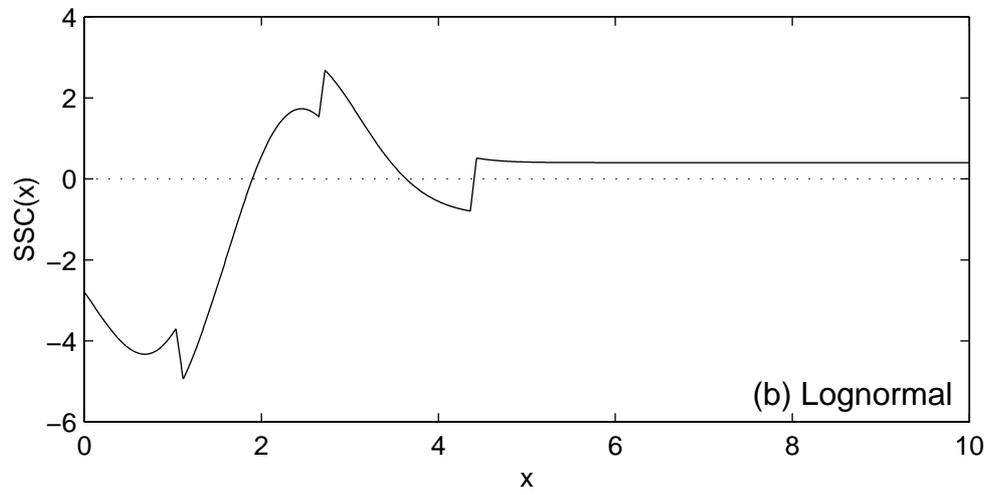

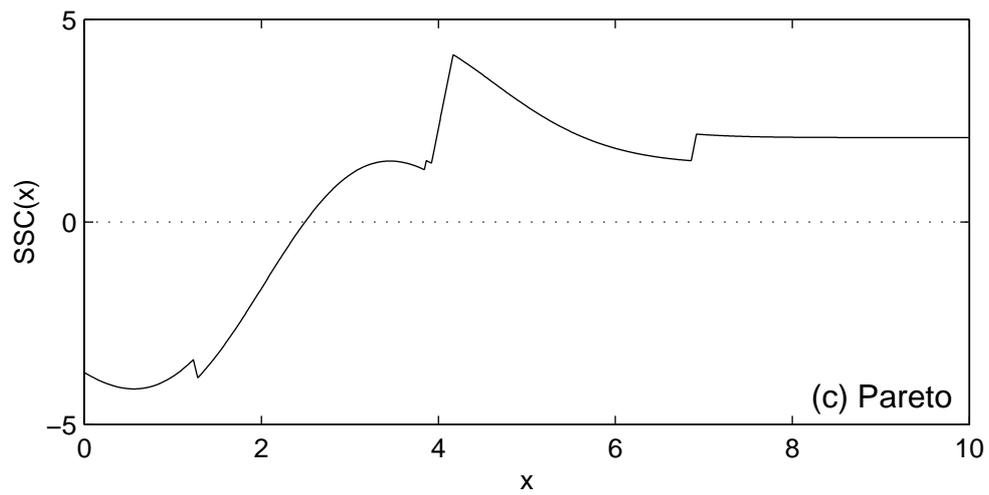

Figure 8

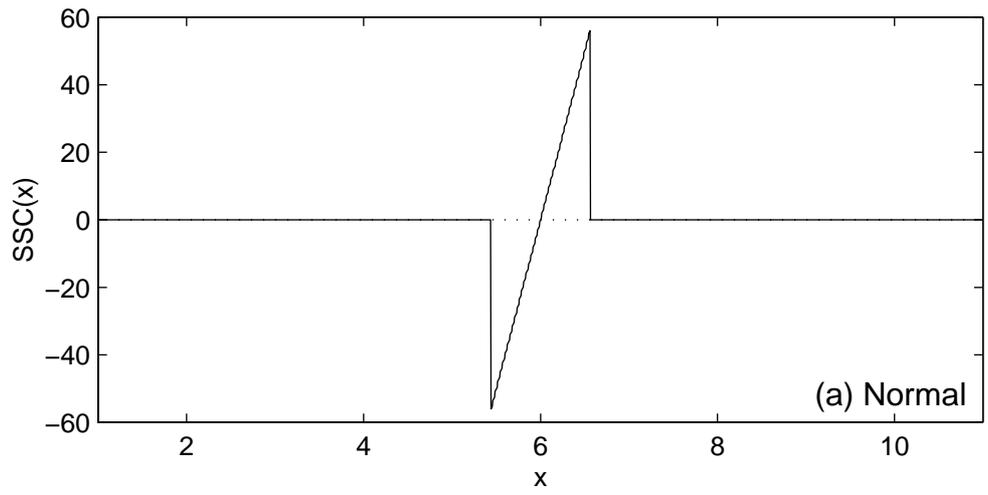
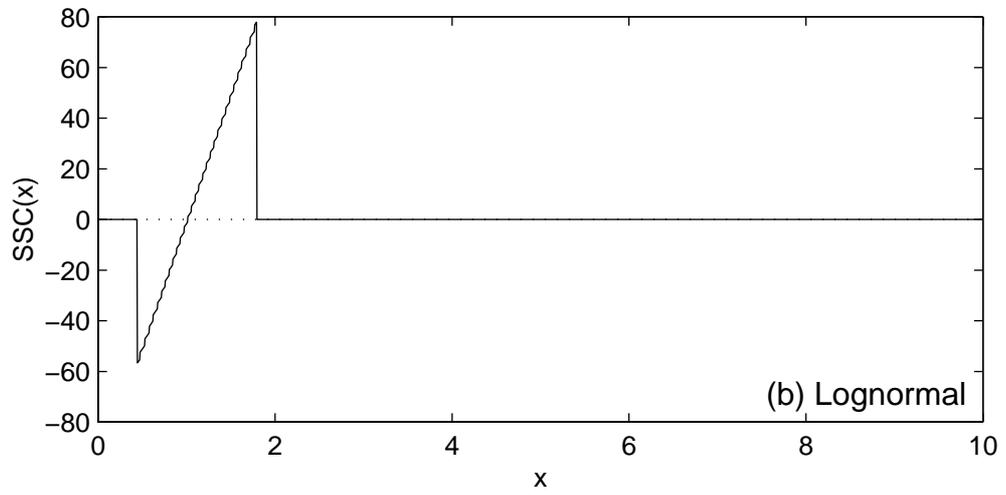
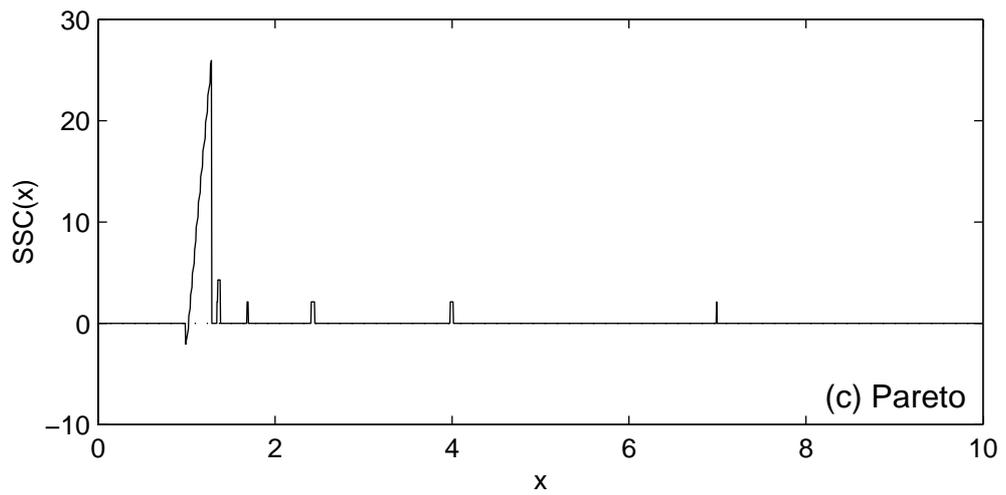

Figure 9

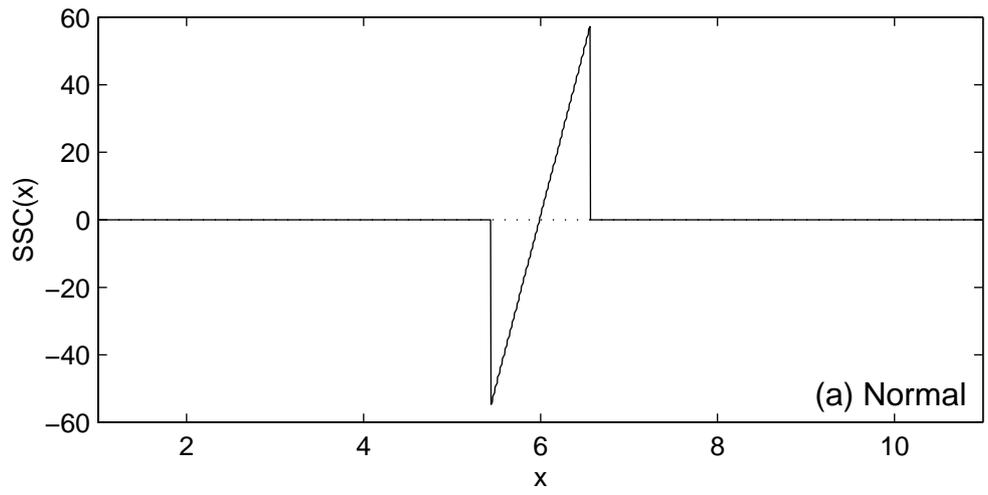

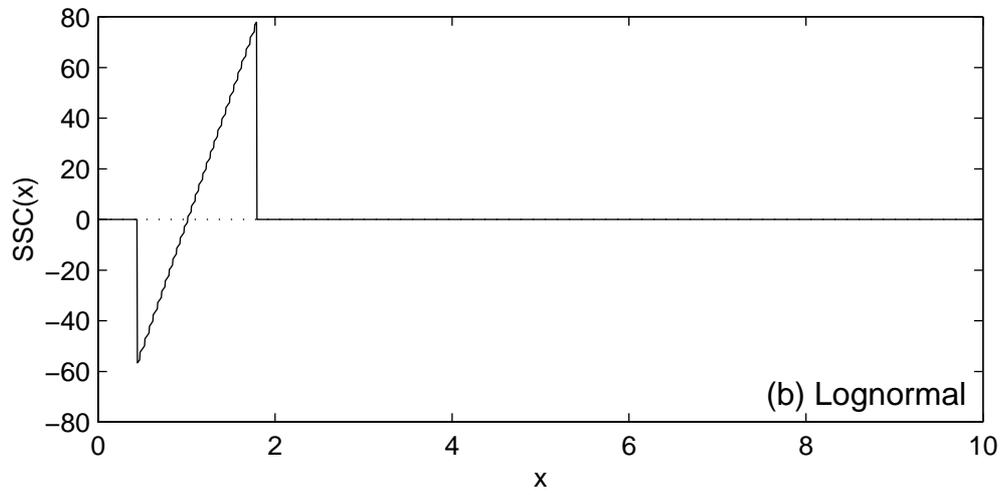

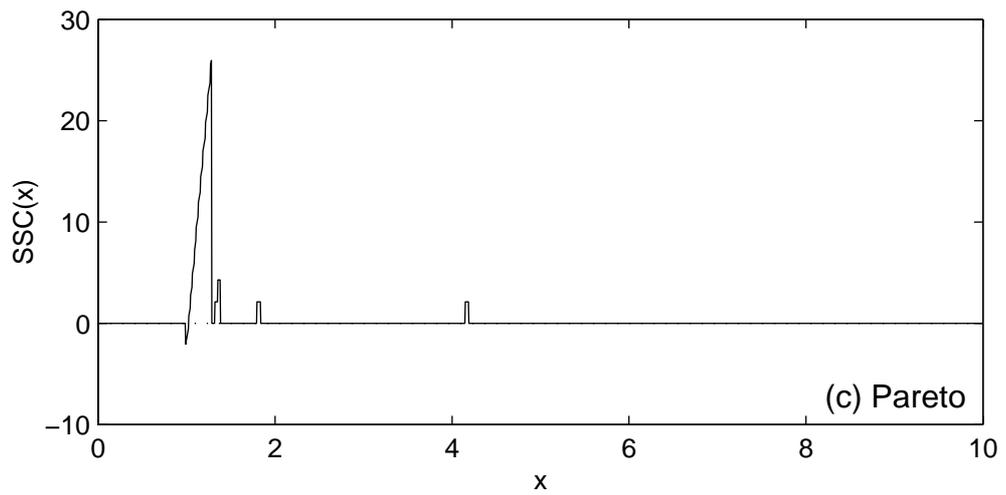

Figure 10

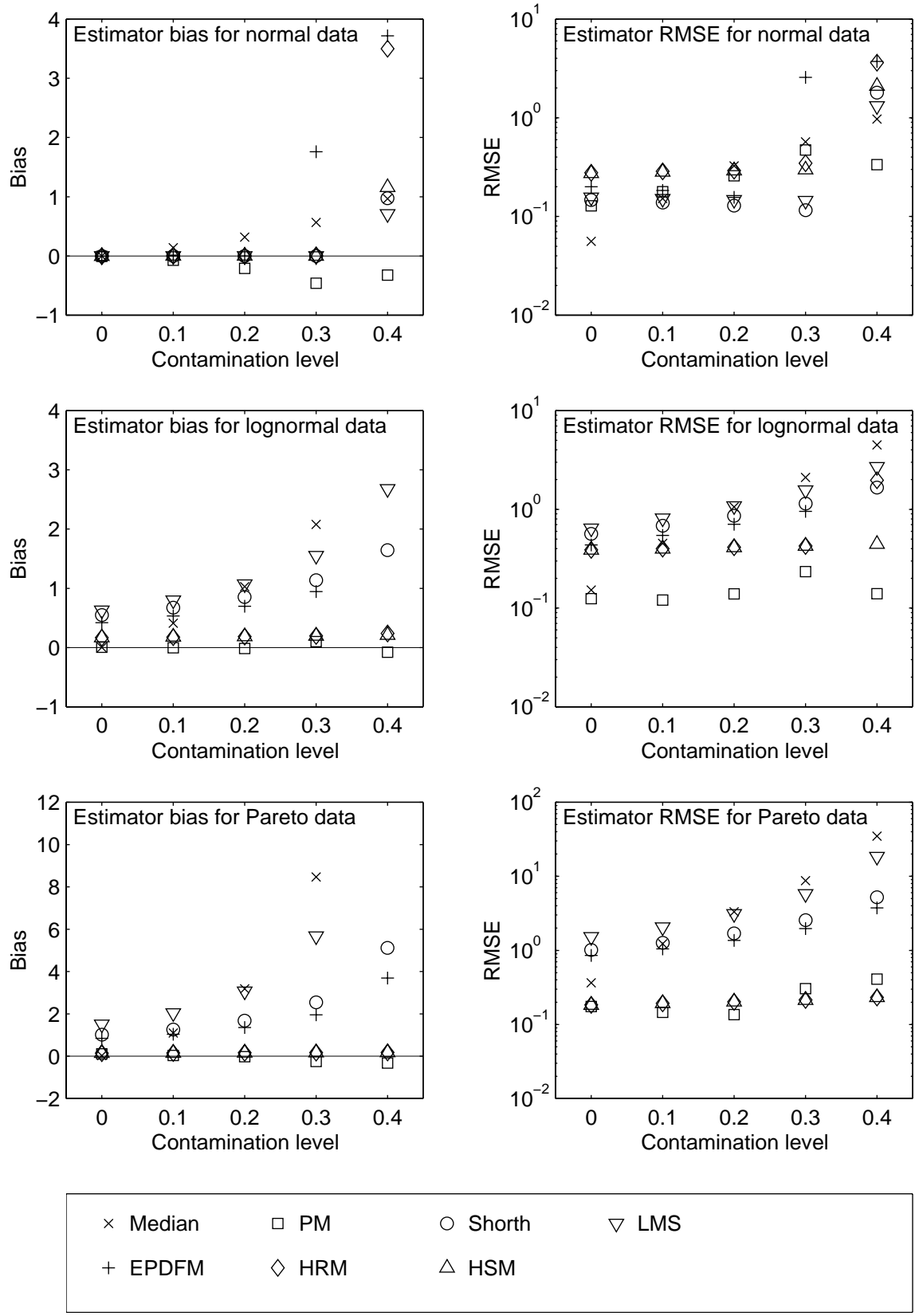

Figure 11

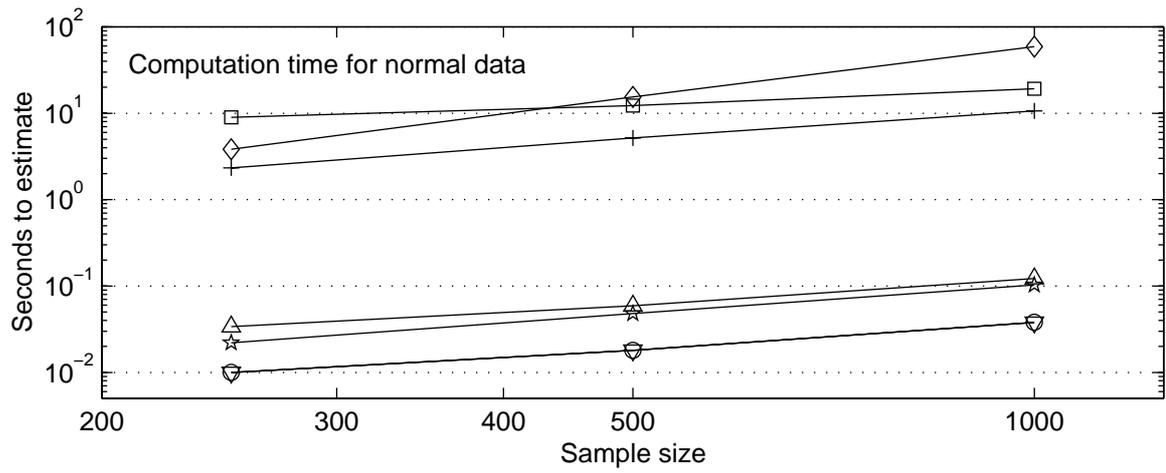

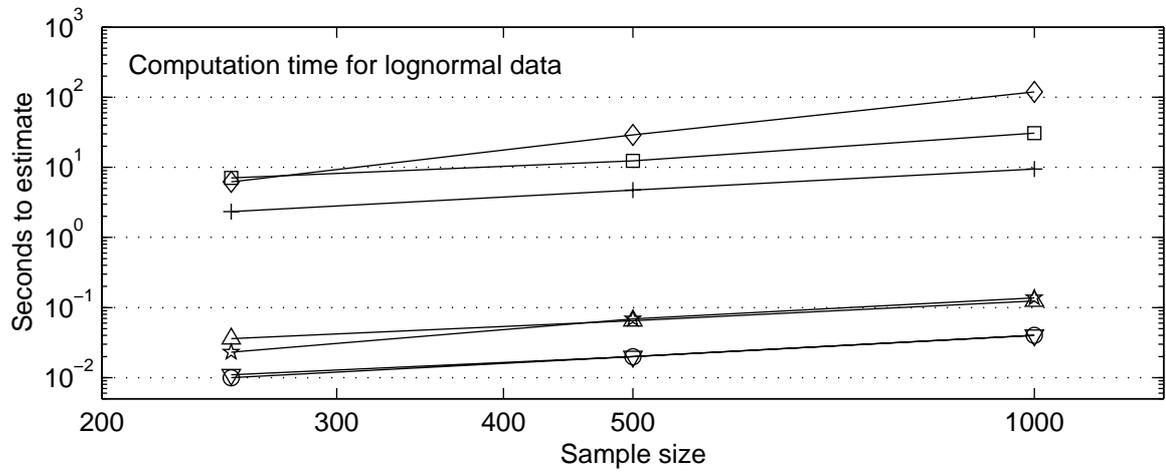

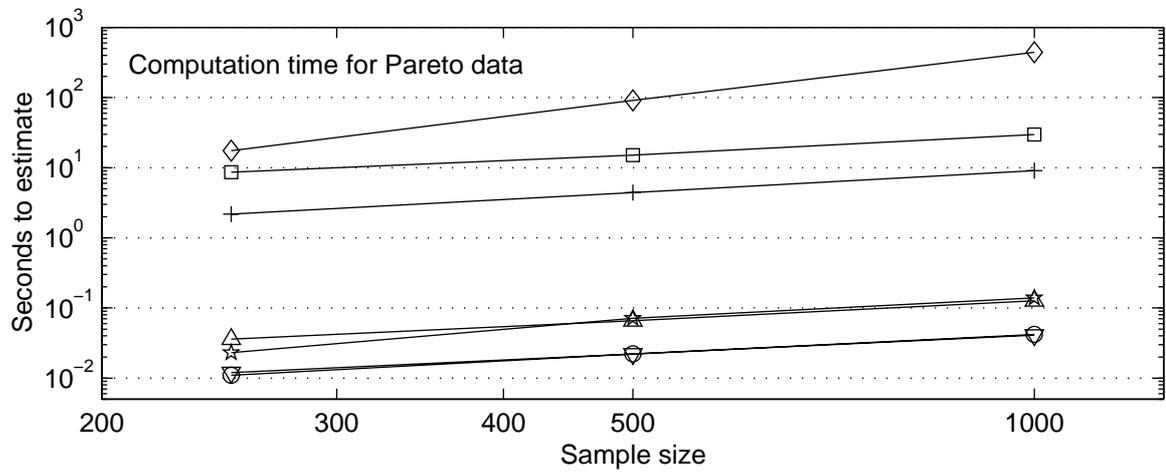

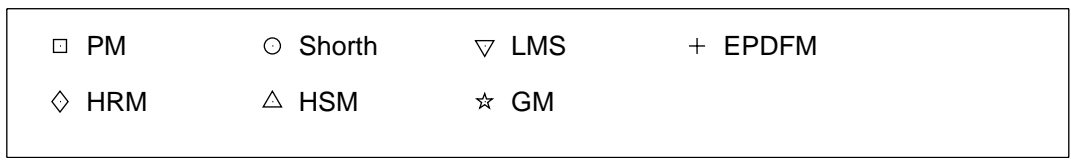

Figure 12

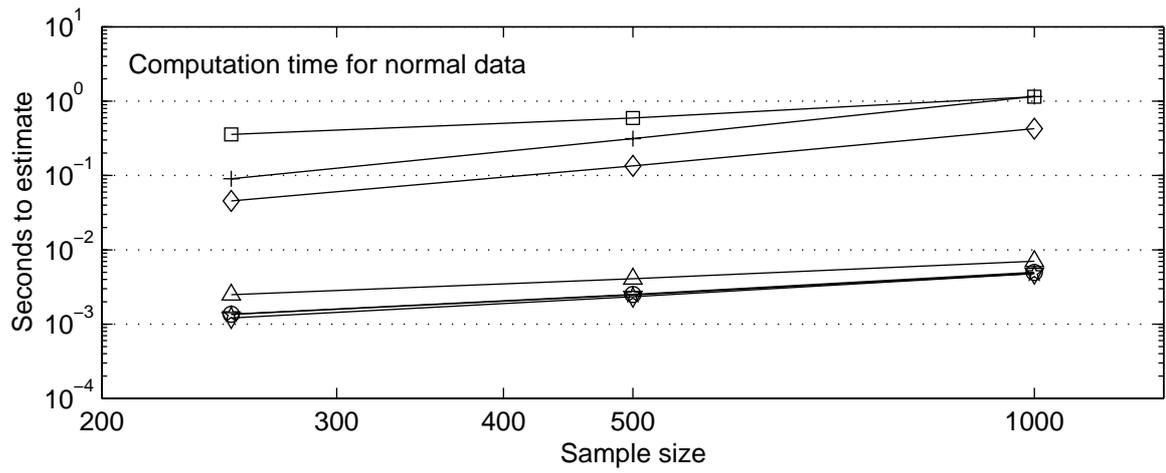

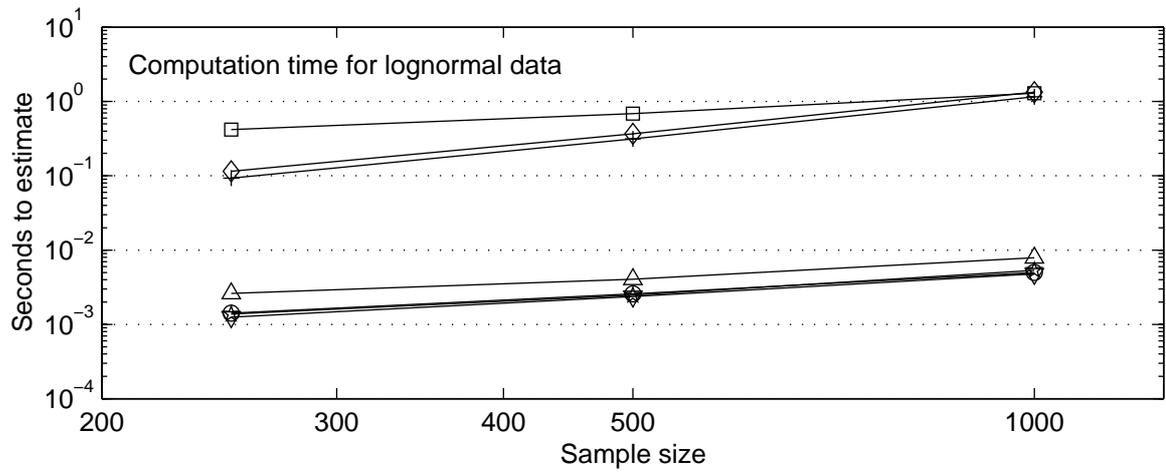

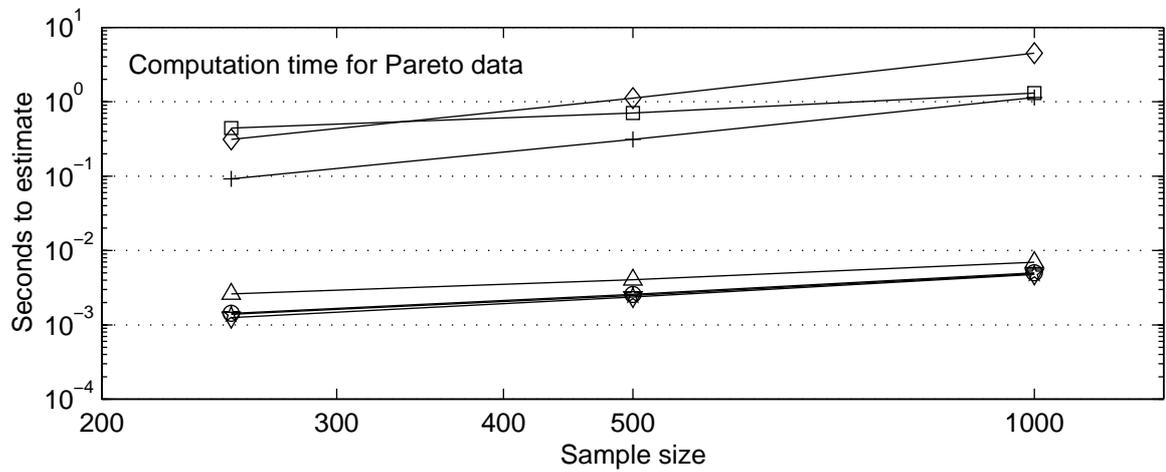

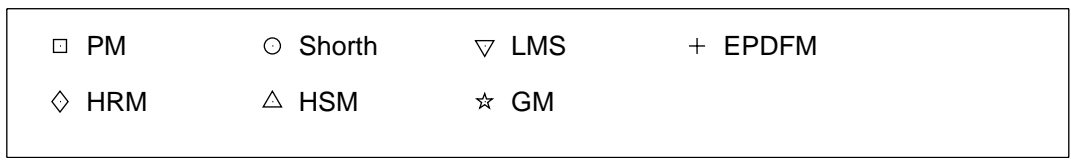

Figure 13

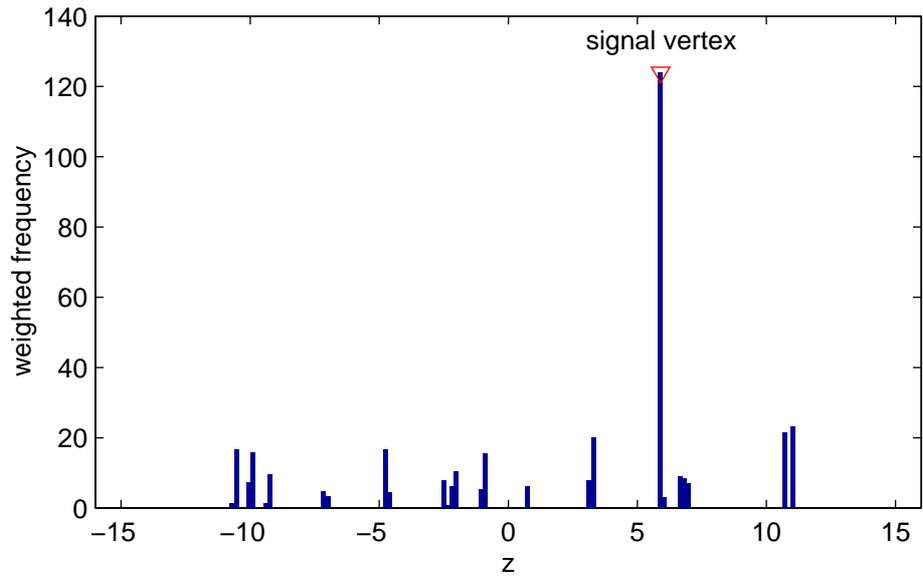

Figure 14

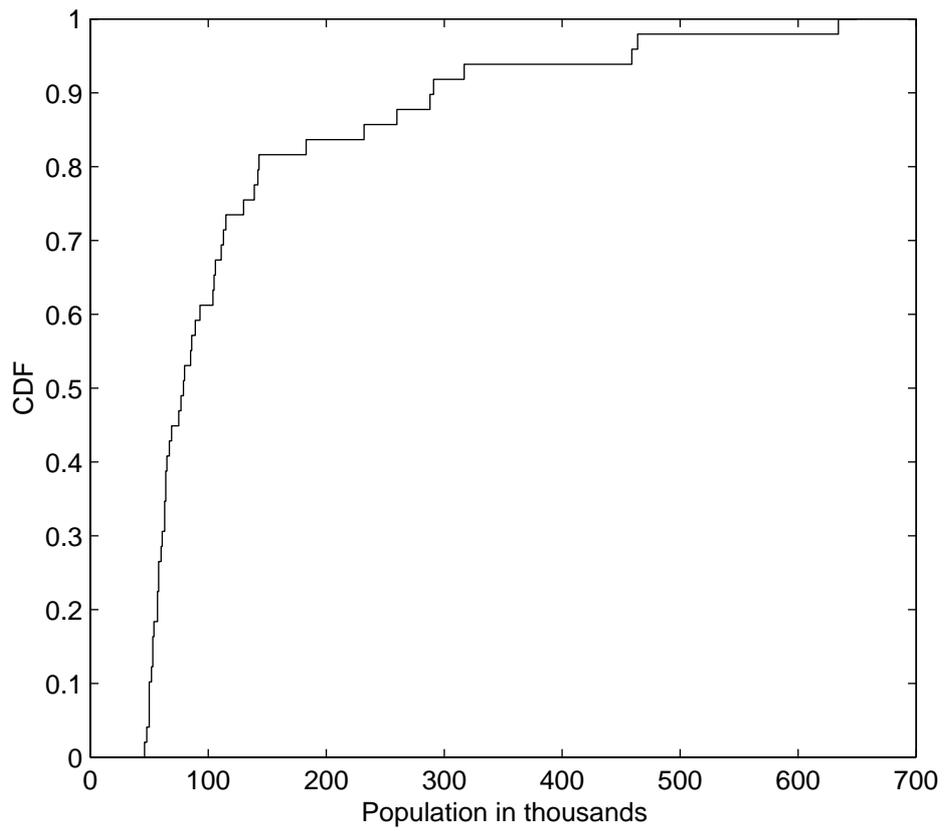

Figure 15

46